\documentclass[a4paper,10pt]{amsart}

\usepackage[T1]{fontenc}
\usepackage[english]{babel}
\usepackage[utf8]{inputenc}
\usepackage{float}
\usepackage{graphicx}
\usepackage{amsmath}
\usepackage{amsfonts}
\usepackage{amssymb}
\usepackage{verbatim}
\usepackage{fullpage}

\usepackage{setspace}

\usepackage{amsthm}

\newtheorem{definition}{Definition}[section]
\newtheorem{theorem}[definition]{Theorem}
\newtheorem{proposition}[definition]{Proposition}

\newtheorem{corollary}[definition]{Corollary}

\author{Emilie Devijver}
\address{Inria Select, Université Paris Sud, Bât. 425, 91405 Orsay Cedex, France}
\email{emilie.devijver@math.u-psud.fr}

\title{Model-based regression clustering for high-dimensional data. Application to functional data.}
\date{\today}
\begin{document}
\maketitle
\begin{abstract}
Finite mixture regression models are useful for modeling the relationship between response and predictors arising from different subpopulations.
In this article, we study high-dimensional predictors and high-dimensional response and propose two procedures to cluster those observations according to the link between predictors and the response.
To reduce the dimension, we propose to use the Lasso estimator, which takes into account the sparsity and a maximum likelihood estimator penalized by the rank, to take into account the matrix structure.
To choose the number of components and the sparsity level, we construct a model collection, varying those two parameters and select a model among this collection with a non-asymptotic criterion.
We extend these procedures to functional data, where predictors and responses are functions. For this purpose, we use a wavelet-based approach.
For each situation, we provide algorithms and apply and evaluate our methods both on simulated and real dataset, to understand how it works in practice.
\keywords{
Model-based clustering, regression, high-dimension, functional data.}
\end{abstract}

\section{Introduction}
Owing to the increasing availability of high-dimensional datasets, regression models for multivariate response and high-dimensional predictors have become important tools.

In this article, we describe two procedures where a random target variable $Y \in \mathbb{R}^q$ depends on explanatory variables or covariates within a cluster-specific regression model.
Each cluster is represented by a parametric conditional distribution, the entire dataset being modeled by a mixture of these distributions. 
It provides a rigorous statistical framework.
The model assumes that each observation $i \in \{1,\ldots,n\}$ originates from one of $K$ disjoint classes and that the data $(Y_i,X_i)$ are independent and identically distributed such that if $i$ belongs to class $k \in \{1,\ldots,k\}$, the
target variable $Y_i$ results from a conditional regression model 
\begin{equation}
\label{modeleFMR}
\boldsymbol{Y}_i=\mathbf{B}_k \boldsymbol{X}_i + \boldsymbol{\varepsilon}_i
\end{equation}
with an unknown matrix of coefficients $\mathbf{B}_k$ and independent errors $\boldsymbol{\varepsilon}_i \in \mathcal{N}_q(0,\Sigma_k)$ with an unknown diagonal covariance matrix $\Sigma_k$.

We work with high-dimensional dataset, in other words the number of parameters to estimate  $ K(qp+q+1)$ is larger than the number of observed target values $q \times n$.
Two ways are considered in this paper, coefficients sparsity and ranks sparsity. 

The first approach consists in estimating the matrix $\mathbf{B}_k$ by a matrix with few nonzero coefficients. 
The well-known Lasso estimator, introduced by \cite{Tibshirani} for linear models, is the solution chosen here.
We refer to  \cite{VdGbook} for an overview of the Lasso estimator and to  \cite{Meinshausen2} for stability selection results.

In the second approach, we consider the  rank sparsity in $\mathbf{B}_k$.
This approach dates back to the $1950$'s and was 
initiated by \cite{anderson} for the linear model. \cite{Izenman} introduced the term of reduced-rank regression for this class of models.
For more recent works, we refer to \cite{Giraud} and to \cite{Bunea1}.
Nevertheless, the linear model is appropriate for homogeneous observations, which is not always the case.

We extend in this paper those methods to mixture regression models.

An important example of high-dimensional dataset is functional dataset.
We refer to Ramsay and Silverman's book \cite{Ramsay} for an overview.
Moreover, a lot of recent works have been done on regression models for functional datasets:  for example, we refer to \cite{Ciarleglio} for a study with scalar response and functional regressors.
In this paper, we focus on the projection of functions into a well-suited wavelet basis.
Indeed, they handle many types of functional data, because they represent global and local attributes of functions and can deal for example with discontinuities.
Moreover, a large class of functions can be well represented with few nonzero coefficients for a suitable wavelet, which leads to sparse matrix of coefficients and then to sparse regression matrix.

We propose here two procedures for clustering high-dimensional or functional datasets, where the high-dimension or functional random target variable $Y \in \mathbb{R}^q$ depends on high-dimensional 
or functional predictors $X \in \mathbb{R}^p$ with a cluster-specific regression model.

Remark that we estimate the number of components, model parameters and cluster proportions. 
In the case of a large number of regressor variables, we use variable selection tools in order to detect relevant regressors.
Since the structure of interest may often be contained into a subset of available variables and many attributes may be useless or even harmful to detect a reasonable clustering structure,
it is important to select the relevant variables.
Moreover, removing irrelevant variables enables us to get an easier  model and can largely enhance comprehension.

Our two procedures are mainly based on three recent works. Firstly, we refer to the article of  \cite{VandeGeer}, which studies finite mixture regression model.
Even if we work on a multivariate version of it, the model considered in the article of \cite{VandeGeer} is adopted here.
The second, the article of \cite{MaugisMeynet}, deals with model-based clustering in density estimation. They propose a procedure, called Lasso-MLE procedure, which determines the number of clusters, the set of relevant variables for the clustering and a clustering of the observations, with high-dimensional data.
We extend this procedure to regression models.
Finally, \cite{Giraud} suggests a low rank estimator for the linear model.
To consider the matrix structure, we develop this last approach to mixture models.

We consider a finite mixture of Gaussian regression models.
The two procedures we propose follow the same sequence of steps.
Firstly, a penalized likelihood approach is considered to determine potential sets of relevant variables. Introduced by \cite{Tibshirani}, the Lasso is used to select variables.
Varying the regularization parameter, it constructs efficiently a data-driven model collection where each model has a reasonable complexity.
The second step of the procedures consists in refitting parameters by a less biased estimator, focusing on selected variables.
Then, we select a model among the collection using the slope heuristic, which was developed by \cite{BirgeMassart}.
The difference between the two procedures is the refitting step. In the first procedure,  later called Lasso-MLE procedure, the maximum likelihood estimator is used/
The second procedure, called Lasso-Rank procedure, deals with low rank estimation. For each model in the collection, a subcollection of models with conditional means estimated by various low rank matrices is constructed. It leads to sparsity for the coefficients and for the rank,
and considers the conditional mean within its matrix structure.

The article is organized as follows. Section \ref{GaussianMixtureRegression} deals with Gaussian mixture regression models. The model collection that we consider is described.
In Section \ref{TwoProcedures}, the two procedures that we propose to solve the problem of high-dimensional regression data clustering are described.
Section \ref{IllustrativeExample} presents an illustrative example, to highlight each choice done in the two procedures.
Section \ref{FunctionalDatasets} states the functional data case, with a description of the projection proposed to convert these functions into coefficients data.
We end this section by studying simulated and benchmark data.
Finally, a conclusion section ends this article.

\section{Finite mixture of Gaussian regression models}
\label{GaussianMixtureRegression}
The model used is a finite Gaussian mixture regression models. 
\cite{VandeGeer} describe this model when the deterministic predictors $\boldsymbol{x} \in \mathbb{R}^p$ are multivariate and the target variable $\boldsymbol{Y} \in \mathbb{R}$ is scalar. 
In this section, it is generalized to multivariate response.
Let us mention that this model has already been introduced, we refer for example to \cite{JonesMcLachlan}. In this paper, we  deal with high-dimensional data.
\subsection{The model}
\label{model}
\label{clustering}
We observe $n$ independent couples $((\boldsymbol{y}_{1}, \boldsymbol{x}_{1}), \ldots, (\boldsymbol{y}_{n}, \boldsymbol{x}_{n}))$,
which are realizations of random variables $((\boldsymbol{Y}_{1}, \boldsymbol{X}_{1}), \ldots, (\boldsymbol{Y}_{n}, \boldsymbol{X}_{n}))$.
Here, $X_i \in \mathbb{R}^p$ are fixed or random covariates and $Y_i \in \mathbb{R}^q$ is a multivariate response variable.
Conditionally to $X_i$,  $Y_i$ is a random variable of unknown conditional density  $s^*(.|X_i)$.
The random response variable $Y_i \in \mathbb{R}^q$ depends on a set of explanatory variables, written $X_i \in \mathbb{R}^p$, through a regression-type model.
If $Y_i$, conditionally to $X_i$, originates from an individual in class $k$, we assume that 

$$\boldsymbol{Y}_i = \mathbf{B}_k \boldsymbol{X}_i + \boldsymbol{\varepsilon}_i ;$$
where $\boldsymbol{\varepsilon}_i \sim \mathcal{N}_q(0, \boldsymbol{\Sigma}_k)$, $\mathbf{B}_k \in M_{p,q}(\mathbb{R})$ is the matrix of class-specific regression coefficients and $\boldsymbol{\Sigma}_k$
is diagonal, positive definite in $\mathbb{R}^q$.

We consider in the following the Gaussian regression mixture model:
\begin{itemize}
 \item Conditionally to the $\boldsymbol{X}_i$'s, $Y_i$'s  are independent, for all $i \in \{1,\ldots, n\}$;
 \item each variable $Y_i$ follows a law of density $s_{\boldsymbol{\xi}}^K(.|\boldsymbol{x}_i)dy$, with
  \begin{align*}
  &s_{\boldsymbol{\xi}}^K(\boldsymbol y|\boldsymbol x)=\sum_{k=1}^{K} \frac{\pi_{k}}{(2 \pi)^{q/2}  \text{det}(\boldsymbol{\Sigma}_k)^{1/2}} \exp \left( -\frac{(\boldsymbol y-\mathbf{B}_{k} \boldsymbol x)^t \boldsymbol{\Sigma}_{k}^{-1}(\boldsymbol y-\mathbf{B}_{k} \boldsymbol x)}{2} \right);\\
  &\boldsymbol{\xi}=(\pi_{1},\ldots, \pi_{K},\mathbf{B}_{1},\ldots,\mathbf{B}_{K},\boldsymbol{\Sigma}_{1},\ldots,\boldsymbol{\Sigma}_{K}) \in \Xi_{K} =  \left( \Pi_{K} \times (\mathbb{R}^{q\times p})^K \times (\mathbb{D}_q^{++})^K \right); \nonumber\\
  & \Pi_{K} = \left\{ (\pi_{1}, \ldots, \pi_{k}) ; \pi_{k} >0 \text{ for } k \in \{1,\ldots, K\} \text{ and } \sum_{k=1}^{K} \pi_{k} = 1 \right\}; \\
  & \mathbb{D}_q^{++} \text{ is the set of diagonal and positive definite matrices in } \mathbb{R}^q .
 \end{align*}
\end{itemize}
We have denoted by $\pi_k$ the proportion of the class $k$.
For all $k \in \{1,\ldots,K\}$, for all $m \in \{1,\ldots,q\}$, for $\boldsymbol{x} = (\xi_1,\ldots, \xi_p) \in \mathbb{R}^p$, 
$[\mathbf{B}_k \boldsymbol{x}]_m= \sum_{j=1}^{p} [\mathbf{B}_{k}]_{m,j} \xi_{j}$ is the $m$th component of the conditional mean of the $k$th mixture component.

We prefer to work with a reparametrized version of this model whose penalized maximum likelihood estimator is scale-invariant and easier to compute.
Indeed, it is not equivariant under scaling of the response.
More precisely, consider the transformation
$$\tilde{Y} = bY, \tilde{\beta}=b\beta, \tilde{\Sigma} = b \Sigma,$$
for $b>0$ which leaves the model invariant.
A reasonable estimator based on transformed data $\tilde{Y}$ should lead to estimator which are related to the first ones up to the homothetic transformation.
This is not the case for the maximum likelihood estimator of $\beta$ and $\Sigma$.
Secondly, the optimization of the loglikelihood is non-convex and hence, it leads to computational issues.
Then, we reparametrize the model described above by generalizing the reparametrization described by St\"adler et al. \cite{VandeGeer}.

For all $k \in \{1,\ldots,K\}$, define new parameters 
\begin{align*}
\boldsymbol{P}_k ^t\boldsymbol{P}_k = \boldsymbol{\Sigma}_k^{-1}\\
\boldsymbol{\Phi}_k=\boldsymbol{P}_k \mathbf{B}_k
\end{align*}
where $\boldsymbol{\Phi}_k \in M_{p,q}(\mathbb{R})$, $P_k$ is the Cholesky decomposition of the positive definite matrix $\boldsymbol{\Sigma}_k^{-1}$. 
Remark that $P_k$ is a diagonal matrix of size $q$.

The model with its reparametrized form then equals

\begin{itemize}
 \item Conditionally to the $\boldsymbol{X}_i$'s, $Y_i$'s  are independent, for all $i \in \{1,\ldots, n\}$;
 \item each variable $Y_i$ follows a law of density $h_\theta^K(.|\boldsymbol{x}_i)dy$, with
 \begin{align*}
  &h_\theta^K(\boldsymbol{y}|\boldsymbol{x})=\sum_{k=1}^{K} \frac{\pi_{k} \det(\boldsymbol{P}_k)}{(2 \pi)^{q/2}} \exp \left( -\frac{(\boldsymbol{P}_k\boldsymbol{y}-\boldsymbol{\Phi}_k \boldsymbol{x} )^t(\boldsymbol{P}_k\boldsymbol{y}-\boldsymbol{\Phi}_k \boldsymbol{x} )}{2} \right)\\
  &\theta=(\pi_{1},\ldots, \pi_{K},\boldsymbol{\Phi}_{1},\ldots,\boldsymbol{\Phi}_{K},\boldsymbol{P}_{1},\ldots,\boldsymbol{P}_{K}) \in \Theta_K = \left( \Pi_{K} \times (\mathbb{R}^{p \times q})^K\times (\mathbb{D}_q)^K \right)\nonumber\\ 
  & \Pi_{K} = \left\{ (\pi_{1}, \ldots, \pi_{K}) ; \pi_{k} >0 \text{ for } k \in \{1,\ldots, K\} \text{ and } \sum_{k=1}^{K} \pi_{k} = 1 \right\} \\
  & \mathbb{D}_q^{++} \text{ is the set of diagonal and positive definite matrices in } \mathbb{R}^q .
 \end{align*}
 \end{itemize}
 
For the sample $((\boldsymbol{x}_1,\boldsymbol{y}_1),\ldots, (\boldsymbol{x}_n,\boldsymbol{y}_n))$, the log-likelihood of this model is equal to, 
$$l(\theta,\textbf{x},\textbf{y})=\sum_{i=1}^n \log \left(\sum_{k=1}^{K} \frac{\pi_{k} \det(\boldsymbol{P}_k)}{(2 \pi)^{q/2}} \exp \left( -\frac{(\boldsymbol{P}_k\boldsymbol{y}_i-\boldsymbol{\Phi}_k\boldsymbol{x}_i)^t(\boldsymbol{P}_k\boldsymbol{y}_i-\boldsymbol{\Phi}_k\boldsymbol{x}_i)}{2} \right) \right) ;$$
and the maximum log-likelihood estimator (MLE) is
$$\hat{\theta}^{MLE} := \underset{\theta \in \Theta_K}{\operatorname{argmin}} \left\{ -\frac{1}{n} l (\theta, \mathbf{x},\mathbf{y}) \right\}.$$
 
Since we deal with the high-dimension case ($p>>n$), this estimator has to be regularized to get stable estimates.
Therefore, we propose the $\ell_1$-norm penalized MLE

\begin{equation}
 \hat{\theta}^{\text{Lasso}}(\lambda) := \underset{\theta \in \Theta_K}{\operatorname{argmin}} \left\{ -\frac{1}{n} l_{\lambda}(\theta, \mathbf{x},\mathbf{y}) \right\} ;
 \label{lassomultidim}
 \end{equation}
where
$$- \frac{1}{n} l_\lambda (\theta, \mathbf{x},\mathbf{y}) = -\frac{1}{n} l(\theta, \mathbf{x},\mathbf{y}) + \lambda \sum_{k=1}^{K} \pi_k ||\boldsymbol{\Phi}_k||_1 ;$$
where $||\boldsymbol{\Phi}_k||_1 = \sum_{j=1}^p \sum_{m=1}^q |[\boldsymbol{\Phi}_k]_{m,j}|$ and with $\lambda > 0$ to specify. This estimator is not the usual $\ell_1$-estimator,
called the Lasso estimator, introduced by \cite{Tibshirani}. It penalizes 
the $\ell_1$-norm of the coefficients matrices $\boldsymbol{\Phi}_k$ and small variances simultaneously, which has some close relations to the Bayesian Lasso (see \cite{ParkCasella}).
Moreover, the reparametrization allows us to consider non-standardized data.

Notice that we restrict ourselves in this article to diagonal covariance matrices which are dependent of the clusters.
We assume that the coordinates of the $Y_i$ are independent, which is a strong assumption, but it allows us to reduce easily the dimension.
We refer to \cite{CeleuxGovaert} for different parametrization of the covariance matrix. Nevertheless, because we work with high-dimensional data ($p$ and $q$ are high and $n$ is small) we prefer a parsimonious model from the diagonal family. We allow different volume clusters because they are capable to detect many clustering structures, as explained in \cite{CeleuxGovaert}. 
 
 \subsection{Clustering with finite mixture of Gaussian regression models}
Suppose that there is a known number of clusters $k$ and assume that we get, from the observations,  an estimator $\hat{\theta}$ 
such that $h_{\hat{\theta}}^K$ well approximates the 
unknown conditional density $s^{*}(.|X_i)$.
We look at this problem as a missing data problem: if we denote by $\boldsymbol{Z}=(\boldsymbol{Z}_{1}, \ldots, \boldsymbol{Z}_{n})$ the component membership,
with
$\boldsymbol{Z}_{i}=([\boldsymbol{Z}_{i}]_1,\ldots, [\boldsymbol{Z}_{i}]_K)$ for $i \in \{1,\ldots,n\}$ is defined by

 $$ [\boldsymbol{Z}_{i}]_k= \left\{ \begin{array}{ll}
                     1 & \text{ if } i \text{ originates from mixture component } k; \\
                     0 & \text{ otherwise;}
                    \end{array}
                    \right.  $$
the complete data are   $((\boldsymbol{x}_{1},\boldsymbol{y}_{1},\boldsymbol{z}_{1}), \ldots, (\boldsymbol{x}_{n}, \boldsymbol{y}_{n},\boldsymbol{z}_{n}))$.

Thanks to the estimation $\hat{\theta}$, we use the Maximum A Posteriori principle (MAP principle) to cluster data.
Specifically, for all $i \in \{1,\ldots,n\}$, for all $k \in \{1,\ldots, K\}$, consider

$$\hat{\boldsymbol{\tau}}_{i,k} (\hat{\theta}) =\frac{\hat{\pi}_{k} \det(\hat{\boldsymbol{P}}_k) \exp \left(- \frac{1}{2} (\hat{\boldsymbol{P}}_k \boldsymbol{y}_i - \hat{\boldsymbol{\Phi}}_k\boldsymbol{x}_i)^t (\hat{\boldsymbol{P}}_k \boldsymbol{y}_i - \hat{\boldsymbol{\Phi}}_k\boldsymbol{x}_i)  \right)}
{\sum_{r=1}^{K} \hat{\pi}_r  \det(\hat{\boldsymbol{P}}_r) \exp \left(- \frac{1}{2} ( \hat{\boldsymbol{P}}_r \boldsymbol{y}_i - \hat{\boldsymbol{\Phi}}_r\boldsymbol{x}_i)^t  ( \hat{\boldsymbol{P}}_r \boldsymbol{y}_i - \hat{\boldsymbol{\Phi}}_r\boldsymbol{x}_i)\right)}$$
the posterior probability of $\boldsymbol{y}_i$ with $i$ coming from the component $k$, where $\theta = (\boldsymbol{\pi},\boldsymbol{\Phi},\textbf{P})$. Then,  data are partitioned by the following rule:

 $$ [\hat{\boldsymbol{Z}}_{i}]_k= \left\{ \begin{array}{ll}
                    1 & \text{ if } \hat{\boldsymbol{\tau}}_{i,k}(\hat{\theta}) > \hat{\boldsymbol{\tau}}_{i,l}(\hat{\theta}) \text{ for all } l\neq k \text{ ;}\\
                    0 & \text{ otherwise.}
                    \end{array}
                    \right. $$

\subsection{Numerical optimization}
\label{EM}
First, we introduce a generalized EM algorithm to approximate $ \hat{\theta}^{\text{Lasso}}(\lambda)$ defined in \eqref{lassomultidim}. 
Then, we discuss initialization and stopping rules in practice, before studying convergence of this algorithm.
\subsubsection{GEM algorithm for approximating the Lasso estimator}
\label{algoGEM}
From an algorithmic point of view, we use a generalized EM algorithm to compute the MLE and the $\ell_1$-norm penalized MLE.
The EM algorithm was introduced by \cite{EM} to approximate the maximum likelihood estimator of mixture model parameters.
 It is an iterative process based on the minimization of the expectation of the likelihood for the complete data conditionally 
 to the observations and to the current estimation of the parameter $\theta^{(\text{ite})}$ at each iteration $\text{(ite)} \in \mathbb{N}^{*}$.
Using to the Karush-Kuhn-Tucker conditions, we  extend the second step to compute the maximum likelihood estimators, penalized or not, under rank constraint or not, as it was done in the scalar case in \cite{VandeGeer}.
All those computations are available in Appendix \ref{em}.
We therefore obtain the next updating formulae for the Lasso estimator defined by \eqref{lassomultidim}. 
Remark that it includes maximum likelihood estimator.
Let denote by $<.,.>$ the euclidean ineer product.
If we denote, for all $j \in \{1,\ldots,p\}$, for all $k \in \{1,\ldots, K\}$, for all $m\in \{1,\ldots,q\}$, for all $i\in \{1,\ldots,n\}$,
\begin{align}
&[\boldsymbol{\widetilde{y}}_i]^{\text{(ite)}}_{k,m} = \sqrt{\boldsymbol{\tau}^{(\text{ite})}_{i,k}} [\boldsymbol{y}_i]_m;\\\nonumber
&[\boldsymbol{\widetilde{x}}_i]^{\text{(ite)}}_{k,j} = \sqrt{\boldsymbol{\tau}^{(\text{ite})}_{i,k}} [\boldsymbol{x}_i]_j ;\\\nonumber
&\Delta_{k,m} = \left( -n_k \langle [\mathbf{\boldsymbol{\widetilde{y}}}_.]^{\text{(ite)}}_{k,m},[\boldsymbol{\Phi}_k]^{\text{(ite)}}_{m,.} [\mathbf{\boldsymbol{\widetilde{x}}}_.]^{\text{(ite)}}_{k,.} \rangle \right)^2
- 4  ||[\mathbf{\boldsymbol{\widetilde{y}}}_.]^{\text{(ite)}}_{k,m}||_ 2^2 ; \\
\label{S}
&[\boldsymbol{S}_{k}]_{j,m}^{(\text{ite})}=-\sum_{i=1}^{n} [\boldsymbol{\widetilde{x}}_i]^{\text{(ite)}}_{k,j} [\boldsymbol{P}_k]^{(\text{ite})}_{m,m} [\boldsymbol{\widetilde{y}}_i]^{\text{(ite)}}_{k,m} +
\sum_{\genfrac{}{}{0pt}{}{j_2=1 }{ j_2\neq j}}^{p} [\boldsymbol{\widetilde{x}}_i]^{\text{(ite)}}_{k,j} [\boldsymbol{\widetilde{x}}_i]^{\text{(ite)}}_{k,j_2} [\boldsymbol{\Phi}_k]^{(\text{ite})}_{m,j_2} ;\\ \nonumber
&n_{k} = \sum_{i=1}^{n} \boldsymbol{\tau}_{i,k}^{(\text{ite})} ; \\ \nonumber
\end{align}
and if we denote by $t^{(ite)} \in (0,1]$ the largest value in $\{0.1^l, l\in \mathbb{N} \}$ such that the update of $\pi$ leads to improve the expected complete penalized log likelihood, 
we  update the parameters by
\begin{eqnarray}
\boldsymbol{\tau}^{(\text{ite})}_{i,k}&=&\frac{\pi_{k}^{(\text{ite})} \left(\det \boldsymbol{P}_k^{(\text{ite})}\right) 
\exp \left( -1/2\left( \boldsymbol{P}^{(\text{ite})}_k \boldsymbol{y}_i-\boldsymbol{\Phi}_k^{(\text{ite})} \boldsymbol{x}_i\right)^t\left(\boldsymbol{P}^{(\text{ite})}_k \boldsymbol{y}_i-\boldsymbol{\Phi}_k^{(\text{ite})} \boldsymbol{x}_i\right) \right)}
{\sum_{r=1}^{K}\pi_{r}^{(\text{ite})}\left(\det \boldsymbol{P}_{r}^{(\text{ite})}\right) \exp \left( -1/2\left(\boldsymbol{P}^{(\text{ite})}_r \boldsymbol{y}_i-\boldsymbol{\Phi}_{r}^{(\text{ite})} \boldsymbol{x}_i\right)^t\left(\boldsymbol{P}^{(\text{ite})}_r \boldsymbol{y}_i-\boldsymbol{\Phi}_{r}^{(\text{ite})}\boldsymbol{x}_i\right) \right)} ; \label{gamma}\\
\label{pi}
\pi_{k}^{\text{(ite+1)}}&=& \pi_{k}^{(\text{ite})}+t^{\text{(ite)}} \left(\frac{n_k}{n} - \pi_{k}^{(\text{ite})}\right) ;\\
\label{rho}
[\boldsymbol{P}_k]_{m_1,m_2}^{\text{(ite+1)}} &=& \left\{ \begin{array}{lll}
     &\frac{n_k \langle [\mathbf{\boldsymbol{\widetilde{y}}}_.]^{\text{(ite)}}_{k,m}, [\boldsymbol{\Phi}_k]_{m,.}^{(\text{ite})} 
[\mathbf{\boldsymbol{\widetilde{x}}}_.]^{\text{(ite)}}_{k,.} \rangle + \sqrt{\Delta}_{k,m}}{2 n_k ||[\mathbf{\boldsymbol{\widetilde{y}}}_.]^{\text{(ite)}}_{k,m}||_2^2}& \text{ if } m_1=m_2=m;\\  
&0 &\text{ elsewhere;}\\
                                                      \end{array}
                                \right.\\
\label{phi}
[\boldsymbol{\Phi}_k]_{m,j}^{\text{(ite+1)}} &=& \left\{ \begin{array}{lll}
            &\frac{-[\boldsymbol{S}_{k}]^{(\text{ite})}_{j,m}+ n\lambda \pi^{(\text{ite})}_{k} }{||[\mathbf{\boldsymbol{\widetilde{x}}}_.]^{\text{(ite)}}_{k,j}||_{2}^2} & \text{ if }  [\boldsymbol{S}_{k}]^{(\text{ite})}_{j,m}>n \lambda \pi^{(\text{ite})}_{k} ;\\
            &-\frac{[\boldsymbol{S}_{k}]^{(\text{ite})}_{j,m}+ n\lambda \pi^{(\text{ite})}_{k} }{||[\mathbf{\boldsymbol{\widetilde{x}}}_.]^{\text{(ite)}}_{k,j}||_{2}^2} & \text{ if } [\boldsymbol{S}_{k}]^{(\text{ite})}_{j,m} < -n \lambda \pi_{k}^{(\text{ite})} ; \\
            & 0 &\text{ else ;}
                                \end{array}
                                \right.
\end{eqnarray}

In our case, the EM algorithm corresponds to alternate between the E-step which corresponds to the computation of \eqref{gamma} and the M-step, which corresponds to the computation of \eqref{pi}, \eqref{rho} and \eqref{phi}.

Remark that to approximate the MLE under rank constraint, we use an easier EM algorithm which is described in details in Appendix.
In the E-step, we compute the a posteriori probability of each observation to belong to each cluster according to the current estimations.
In the M-step, we consider each observation within its estimated cluster, and consider the linear regression estimators under rank constraint by keeping the biggest singular values.
\subsubsection{Initialization and stopping rules}
We initialize the clustering process with the $k$-means algorithm on couples $(\boldsymbol{x}_i,\boldsymbol{y}_i) \in \mathbb{R}^{p+q}$, for all $i \in \{1,\ldots,n\}$.
According to this clustering, we compute the linear regression estimators in each class.
Then, we run a small number of times the EM-algorithm, repeat this initialization many times and keep the one which corresponds to the highest log-likelihood.

To stop the algorithm, we propose to run it a minimum number of times, and to specify a maximum number of iterations to make sure that it will stop.

Between these two bounds, we stop if a relative criteria on the log-likelihood is small enough and if a relative criteria on the parameters is small enough too.
Those criteria are adapted from \cite{VandeGeer}.

By those rules, we are trying to attain a local optimum as close to a global optimum of the likelihood function as possible.

\subsubsection{Numerical convergence of this algorithm}
We address here the  convergence properties of the algorithm described in paragraph \ref{algoGEM}. 
Although convergence to stationary points has been get for the EM algorithm (see \cite{Wu}), it is not true without conditions hard to verify for generalized EM algorithm.

Results we propose are quite similar to the one get by \cite{VandeGeer}, because algorithms are similar.
We refer the interested reader to this article for proof of the following results, when $q=1$. Same ideas are working for any values of $q$.

\begin{proposition}
\label{prop1}
The algorithm described in paragraph \ref{algoGEM} has the descent property:
\begin{align*}
-\frac{1}{n} l_{\lambda} (\theta^{(\text{ite}+1)},\mathbf{x},\mathbf{y}) \leq 
-\frac{1}{n} l_{\lambda} (\theta^{(\text{ite})},\mathbf{x},\mathbf{y})
\end{align*}
\end{proposition}
The Proposition \ref{prop1} is clear by construction of $\theta^{(\text{ite}+1)}$, the M-step of the algorithm as a coordinate-wise minimization.

\begin{proposition}
\label{prop2}
Assume that $Y_i \neq 0$ for all $i \in \{1,\ldots,n\}$. Then, for $\lambda >0$, 
the function $\theta \mapsto l_{\lambda} (\theta,\mathbf{x},\mathbf{y})$ is bounded from above for all values $\theta \in \Theta_K$.
\end{proposition}
This is true because we penalize the log-likelihood by $\boldsymbol{\Phi}=(\boldsymbol{\Phi}_1,\ldots, \boldsymbol{\Phi}_K)$, where $\boldsymbol{\Phi}_k =\boldsymbol{P}_k \mathbf{B}_k$,
then small variances are penalized also.
The penalized criteria therefore stays finite whenever $\Sigma_k$, $k \in \{1,\ldots,K\}$.
This is not the case for the unpenalized, which is unbounded if the variance of a variable tends to $0$. We refer for example to \cite{McLachlanPeel}.

\begin{corollary}
\label{cor}
For the algorithm described in paragraph \ref{algoGEM},
$- \frac{1}{n} l_{\lambda} (\theta^{\text{(ite)}},\mathbf{x},\mathbf{y})$ decreases monotonically to some value $\bar{l}> -\infty$.

\end{corollary}

According to this corollary, we know that the algorithm converges to some finite value, depending on initial values.

\begin{theorem}
\label{thm}
For the algorithm described in paragraph \ref{algoGEM}, every cluster point $\bar{\theta} \in \Theta_K$ of the sequence $\left\{ \theta^{\text{(ite)}}, \text{ite}=0,1,\ldots \right\}$, generated by the algorithm, is a stationary point of the function
$\theta \mapsto l_\lambda(\theta,\mathbf{x},\mathbf{y})$.
\end{theorem}
This theorem proves that the algorithm converges to a stationary point.
We refer to \cite{Tseng} for the definition of a stationary point for non-differentiable functions.

Corollary \ref{cor} is deduced from Proposition \ref{prop1} and Proposition \ref{prop2}.
Theorem \ref{thm} is more difficult to prove.

\subsection{Variable selection and model collection}
\label{varselect}
We deal with high-dimensional data where we observe a sample of small size $n$ and we have to estimate many coefficients ($K(pq+q+1)-1$).
Then, we have to focus on variables that are relevant for the clustering.
The notion of irrelevant indices has to be defined.
\begin{definition}
A couple $(x_j,y_z)$ is said to be \emph{irrelevant} for the clustering if 
$[\boldsymbol{\Phi}_{1}]_{m,j} = \ldots = [\boldsymbol{\Phi}_{K}]_{m,j}=0$, which means that the variable $x_j$ does not explain the variable $y_z$ for the clustering process.
We also say that the indices $(j,z)$ is irrelevant if the couple $(x_j,y_z)$ is irrelevant.
A \emph{relevant} couple is a couple which is not irrelevant: at least in one cluster $k$, the coefficient $[\boldsymbol{\Phi}_k]_{m,j}$ is not equal to zero.
We denote by $J$ the indices set of relevant couples.
\end{definition}

Remark that $J \subset \{1,\ldots,q\}\times \{1,\ldots,p\}$.
We denote by $^c J$ the complement of $J$ in $\{1,\ldots,q\} \times \{ 1,\ldots,p\}$.
For all $k \in \{1,\ldots,K\}$, we denote by $\boldsymbol{\Phi}_k^{[J]}$ the matrix of size $q\times p$ with $0$ on the set $^c J$.
Be also define by 
$\mathcal{H}_{(K,J)}$ the model with $K$ components and with $J$ for indices set of relevant couples:

\begin{align}
\label{modele h}
\mathcal{H}_{(K,J)} &= \left\{ h_{\theta}^{(K,J)}(\boldsymbol{y}|\boldsymbol{x}) = \sum_{k=1}^{K} \frac{\pi_{k} \det(\boldsymbol{P}_k)}{(2 \pi)^{q/2}} \exp \left( -\frac{(\boldsymbol{P}_k\boldsymbol{y}-\boldsymbol{\Phi}_k^{[J]} \boldsymbol{x} )^t(\boldsymbol{P}_k\boldsymbol{y}-\boldsymbol{\Phi}_k^{[J]} \boldsymbol{x} )}{2} \right), \right.\\
& \hspace{1cm}\left. \theta=(\pi_1,\ldots, \pi_K,\boldsymbol{\Phi}_1^{[J]},\ldots, \boldsymbol{\Phi}_K^{[J]}, \boldsymbol{P}_1,\ldots,\boldsymbol{P}_K) \in \Theta_{(K,J)} = \Pi_K \times \left( \mathbb{R}^{q \times p} \right)^K \times \left(\mathbb{R}_+^{q} \right)^K \right\}. \nonumber
 \end{align}

We construct a model collection by varying the number of components $K$ and the indices set of relevant couples $J$.

The subset $J$ is  constructed with the Lasso estimator defined in \eqref{lassomultidim}. Nevertheless, to be consistent with the definition of relevant couples, the Group-Lasso estimator could be preferred,
where coefficients $\{\boldsymbol{\Phi}_{1}]_{m,j},\ldots,\boldsymbol{\Phi}_{K}]_{m,j}\}$ are groupped.
We focus here on the Lasso estimator, but the Group-Lasso approach is described in Appendix. It gives mainly the same results, but 
the Lasso estimator leads to consider also isolated irrelevant couples.

\section{Two procedures}
\label{TwoProcedures}
The goal of our procedures is, given a sample $(\textbf{x},\textbf{y})=((\boldsymbol{x}_1,\boldsymbol{y}_1),\ldots,(\boldsymbol{x}_n,\boldsymbol{y}_n)) \in (\mathbb{R}^p \times \mathbb{R}^q)^n$, to discover the 
relation between the variable $\boldsymbol{Y}$ and the variable  $\boldsymbol{X}$. 
Thus, we have to estimate, according to the representation of $\mathcal{H}_{(K,J)}$, the number of clusters $K$, the relevant variables set $J$, 
and the parameter $\theta$.
To overcome this difficulty, we want to take advantage of the sparsity property of the $\ell_1$-penalization to perform automatic variable selection in clustering  high-dimensional data.
Then, we compute another estimator restricted on relevant variables, which will work better because it is no longer an high-dimensional issue. Thus, 
we avoid shrinkage problems due to the Lasso estimator.
The first procedure takes advantage of the maximum likelihood estimator, whereas the second one takes into account the matrix structure of $\boldsymbol{\Phi}$ with a low rank estimation.
 \subsection{Lasso-MLE procedure}
 \label{LassoMLE}
This procedure is decomposed into three main steps: we construct a model collection, then in each model we compute the  maximum likelihood estimator and finally we select the best one among the model collection.

The first step consists of constructing a model collection $\{\mathcal{H}_{(K,J)}\}_{(K,J) \in \mathcal{M}}$ in which $\mathcal{H}_{(K,J)}$ is defined by equation
\eqref{modele h},
and the model collection is indexed by $ \mathcal{M}= \mathcal{K} \times \mathcal{J}$. 
We denote by $\mathcal{K}\subset \mathbb{N}^*$ the possible number of components.
We assume that  we could bound $\mathcal{K}$ without loss of generality. We also note $\mathcal{J} \subset \mathcal{P} (\{1,\ldots,q\} \times \{1,\ldots,p\})$.

To detect the relevant variables and construct the set $J \in \mathcal{J}$, we penalize the log-likelihood by an $\ell_1$-penalty on the mean parameters proportional
to $||\boldsymbol{\Phi}_k||_1 = \sum_{j=1}^p \sum_{m=1}^q |[\boldsymbol{\Phi}_k]_{m,j}|$.
In $\ell_1$-procedures, the choice of the regularization parameters is often difficult: fixing the number of components $K \in \mathcal{K}$, we propose to construct a 
data-driven grid  $G_K$ of regularization parameters by using the updating formulae of the mixture parameters in the EM algorithm. We can give a formula for $\lambda$, 
the regularization parameter, depending on which coefficients we want to shrink to zero, for all $k \in \{1,\ldots,K\}, j\in \{1,\ldots,p\}, z\in \{1,\ldots,q\}$:
\begin{align*}
 [\boldsymbol{\Phi}_k]_{m,j} = 0 \hspace{0.5 cm}  \Leftrightarrow& \hspace{0.5 cm}  [\boldsymbol{\lambda}_{k}]_{j,m}= \frac{|[\boldsymbol{S}_{k}]_{j,m}|}{n \pi_{k}} ;
\end{align*}
where $[\boldsymbol{S}_{k}]_{j,m}$ is defined by \eqref{S}.
Then, we define the data-driven grid by 
\begin{align}
G_K = \left\{[\boldsymbol{\lambda}_{k}]_{j,m}, k\in \{1,\ldots,K\}, j \in \{1,\ldots,p\}, m \in \{1,\ldots,q\}\right\}.
\label{grid}
\end{align}
We could compute it from maximum likelihood estimations.

Then, for each $\lambda \in G_K$, we could compute the Lasso estimator defined by
$$\hat{\theta}^{\text{Lasso}}(\lambda) = \underset{\theta \in \Theta_{(K,J)}}{ \operatorname{argmin}} \left\{ -\frac{1}{n} \sum_{i=1}^n 
\log (h_{\theta}^{(K,J)}(\boldsymbol{y}_i|\boldsymbol{x}_i)) + \lambda \sum_{k=1}^K \pi_k ||\boldsymbol{\Phi}_k||_1 \right\} .$$ 
For a fixed number of mixture components $K \in \mathcal{K}$ and a regularization parameter $\lambda \in G_K$, we could use an EM algorithm, recalled in Appendix \ref{em}, to approximate this estimator.
Then, for each $K \in \mathcal{K}$ and for each $\lambda \in G_K$, we could construct the relevant variables set $J_\lambda$. We denote by $\mathcal{J}$ the collection of these sets.

The second step consists of approximating the MLE 
$$\hat{h}^{(K,J)}= \underset{h \in \mathcal{H}_{(K,J)}}{ \operatorname{argmin}} \left\{ -\frac{1}{n} \sum_{i=1}^n \log (h(\boldsymbol{y}_i|\boldsymbol{x}_i)) \right\} ; $$ 
using the EM algorithm for each model $(K,J)\in \mathcal{M}$. 

The third step is devoted to model selection.
Rather than selecting the regularization parameter, we select the refitted model.
Instead of using an asymptotic criterion, as BIC or AIC, we use the slope heuristic described in \cite{BirgeMassart}, which is a non-asymptotic criterion for selecting a model among a model collection.
Let us explain briefly how it works.
Firstly, models are grouping according to their dimension $D$, to obtain a model collection $\{\mathcal{H}_{D} \}_{D \in \mathcal{D}} $. The dimension of a model
is the number of parameters estimated in the model.
For each dimension $D$, let $\hat{h}_{D}$ be the estimator maximizing the likelihood among the estimators associated to a model of dimension $D$.
Also, the function $D/n \mapsto 1/n \sum_{i=1}^n \log (\hat{h}_{D}(\boldsymbol{y}_i|\boldsymbol{x}_i))$ has a linear behavior for large dimensions. We estimate the slope, denoted by $\hat{\kappa}$, which will be used to calibrate the penalty.
The minimizer $\hat{D}$ of the penalized criterion $-1/n \sum_{i=1}^n \log (\hat{h}_{D}(\boldsymbol{y}_i|\boldsymbol{x}_i)) + 2 \hat{\kappa} D/n$ is determined, 
and the model selected is $(K_{\hat{D}},J_{\hat{D}})$.
Remark that $D = K(|J|+q+1)-1$.

Note that the model is selected after refitting, which avoids issue of regularization parameter selection.
For an oracle inequality to justify the slope heuristic used here, see \cite{inegOracleLassoMLE}.

\subsection{Lasso-Rank procedure}
Whereas the previous procedure does not take into account the multivariate structure, we propose a second procedure to perform this point.
For each model belonging to the collection $\mathcal{H}_{(K,J)}$, a subcollection is constructed, varying the rank of $\boldsymbol{\Phi}$. Let us describe this procedure.

As in the Lasso-MLE procedure, we first construct a collection of models, thanks to the $\ell_1$-approach.
For $\lambda \geq 0$, we obtain an estimator for $\theta$, denoted by $\hat{\theta}^{\text{Lasso}}(\lambda)$, for each model belonging to the collection.
We could deduce the set of relevant columns, denoted by $J_\lambda$ and this for all $K \in \mathcal{K}$: we deduce $\mathcal{J}$ 
the collection of relevant variables set.

The second step consists to construct a subcollection of models with rank sparsity, denoted by 
$$\{\check{\mathcal{H}}_{(K,J,R)}\}_{(K,J,R) \in \check{\mathcal{M}}}.$$
The model $\check{\mathcal{H}}_{(K,J,R)}$ has $K$ components, the set $J$ for active variables and $R$ is the vector
of the ranks of the matrix of regression coefficients in each group:

\begin{equation}
\label{modele htilde}
\check{\mathcal{H}}_{(K,J,R)}= \left\{\boldsymbol  y \in \mathbb{R}^q | \boldsymbol x \in \mathbb{R}^p \mapsto h_{\theta}^{(K,J,R)}(\boldsymbol y|\boldsymbol x)   \right\}
 \end{equation}
where
\begin{align*}
 h_{\theta}^{(K,J,R)}(\boldsymbol{y}|\boldsymbol{x}) &= \sum_{k=1}^{K} \frac{\pi_{k} \det(\boldsymbol{P}_k)}{(2 \pi)^{q/2}} \exp \left( -\frac{(\boldsymbol{P}_k\boldsymbol{y}-(\boldsymbol{\Phi}_k^{R_k})^{[J]} \boldsymbol{x} )^t(\boldsymbol{P}_k\boldsymbol{y}-(\boldsymbol{\Phi}_k^{R_k})^{[J]} \boldsymbol{x} )}{2} \right) ;\\
\theta&=(\pi_1,\ldots, \pi_K,(\boldsymbol{\Phi}_1^{R_1})^{[J]},\ldots, (\boldsymbol{\Phi}_K^{R_K}))^{[J]}, \boldsymbol{P}_1,\ldots,\boldsymbol{P}_K) \in \Theta_{(K,J,R)}\\
 \Theta_{(K,J,R)} &= \Pi_K \times \Psi_{(K,J,R)} \times \left(\mathbb{R}_+^{q} \right)^K ;\\
\Psi_{(K,J,R)} &= \left\{\left. ((\boldsymbol{\Phi}_1^{R_1})^{[J]},\ldots, (\boldsymbol{\Phi}_K^{R_K})^{[J]}) \in \left(\mathbb{R}^{q \times p} \right)^K \right| \text{Rank}(\boldsymbol{\Phi}_k)=R_k \text{ for all } k \in \{1,\ldots,K\} \right\} ; 
\end{align*}

and $\check{\mathcal{M}}= \mathcal{K} \times \mathcal{J} \times \mathcal{R}$. 
We denote by $\mathcal{K} \subset \mathbb{N}^*$ the possible number of components, $\mathcal{J}$ a collection of subsets of $\{1,\ldots,p\}$ and $\mathcal{R}$ the set of vectors of size $K \in \mathcal{K}$ with rank values for each mean matrix.
We compute the MLE under the constrained ranks thanks to an EM algorithm.
Indeed, we constrain the estimation of $\boldsymbol{\Phi}_k$, for the cluster $k$, to have a rank equals to $R_k$, by keeping only the $R_k$ largest singular values.
More details are given in Section \ref{EM2}.
It leads to an estimator of the mean with row sparsity and low rank for each model. 
As described in the above section, a model is selected using the slope heuristic.
This step is justified theoretically in \cite{inegOracleLassoRank}.

\subsection{Some comments}
Those two procedures have both been implemented in Matlab, with the help of Benjamin Auder and the Matlab code is available.
We run the EM algorithm several times: once to construct the regularization grid, twice for each regularization parameter for the Lasso-MLE procedure (once to approximate the Lasso estimator,
and once to refit it with the maximum likelihood estimator) and more times for the Lasso-Rank procedure (we vary also the ranks vector).
If we look at every regularization parameters in the grid defined in \eqref{grid}, there are $K*p*q$ values and then we compute the EM algorithm $2*K*p*q +1$ times, which could be large.
Even if each EM algorithm is fast (implemented with C), repeat it numerous times could be time-consuming. We propose to the user to select relevant regularization parameters: either a regular subcollection of $G_K$, to get various sparsities, or focus on the large values of regularization parameters, to get sparse solutions.

The model collection constructed by the Lasso-MLE procedure is included in the model collection constructed by the Lasso-Rank procedure:
$$\{ \mathcal{H}_{(K,J)}\}_{(K,J) \in {\mathcal{M}}} \subset \{\check{\mathcal{H}}_{(K,J,R)}\}_{(K,J,R) \in \check{\mathcal{M}}}.$$
Then, the second procedure will work better.
Nevertheless, it is time-consuming to construct all those models, with various rank vectors and maybe it is not necessary, depending on the data.
It is known that the framework of low-rank  arises  in  many  applications,
among which analysis of EEG data decoding \cite{AndersonStolzShamsunder}, neural
response modeling \cite{Brown}.
However, in some data the low-rank structure could be not adapted and then not needed.
\section{Illustrative example}
\label{IllustrativeExample}
We illustrate our procedures on four different simulated dataset, adapted from \cite{VandeGeer}.
Firstly, we present the  models used in these simulations. Then, we validate numerically each step and we finally compare the results of our procedures with others.
Remark that we propose here some examples to illustrate our methods, but not a complete analysis. We highlight some issues which seem important.
Moreover, we do not illustrate the one-component case, focusing on the clustering.
If, on some dataset, we are not convinced by the clustering, we could add to the model collection models with one component, more or less sparse, using the same pattern (computing
the Lasso estimator to get the relevant variables for various regularization parameters and refitting parameters with the maximum likelihood estimator, under constrained ranks or not) and then we select a model among this collection of linear and mixture models.
\subsection{The model}
Let $\textbf{x}$ be a sample of size $n$ distributed according to multivariate standard Gaussian. We consider a mixture of two components.
Besides, we fix the number of active variables to $4$ in each cluster. 
More precisely, the first four variables of $\boldsymbol{Y}$ are explained respectively by the four first variables of $\boldsymbol{X}$.
Fix $\pi = \left(\frac{1}{2},\frac{1}{2}\right)$ and $\boldsymbol{P}_k=I_q$ for all $k \in \{1,2\}$, where $I_q$ denotes the identity matrix of size $q$.

The difficulty of the clustering is partially controlled by the signal-to-noise ratio (denoted SNR). In this context, we could extend the natural idea of the SNR with the following definition, where $\text{Tr}(A)$ denotes the trace of the matrix $A$.
$$\text{SNR}= \frac{\text{Tr}(Var(\boldsymbol{Y}))}{\text{Tr}(Var(\boldsymbol{Y}|\mathbf{B}_k=0 \text{ for all } k \in \{1,\ldots,K\}))}.$$
Remark that it only controls the distance between the signal with or without the noise and not the distance between the both signals.

We compute four different models, varying $n$, the SNR and the distance between the clusters. Details are available in the Table \ref{tableRecap}.

\begin{table}[H]
\centering

\begin{tabular}{|c|c|c|c|c|c|}
 \hline
  \phantom{bb} &  Model 1 & Model 2 & Model 3 & Model 4 & Model 5\\ \hline 
 \hline
 n & 2000 & 100 & 100 & 100 & 50 \\ \hline
 k & 2 &2&2&2&2\\ \hline
 p & 10 & 10 & 10 & 10 & 30 \\ \hline
 q & 10 & 10 & 10 & 10 & 5 \\ \hline
  $\mathbf{B}_{1|J}$& 3 &3 & 3 &5 & 3 \\ \hline
  $\mathbf{B}_{2|J}$&  -2& -2& -2 & 3& -2\\ \hline
  $\sigma$ & 1 & 1 & 3 & 1& 1\\ \hline
  SNR & 3.6 & 3.6 & 1.88 & 7.8&3.6\\ \hline
\end{tabular}
\caption{Description of the different models.}
\label{tableRecap}

\end{table}

Take a sample of $\boldsymbol{Y}$ according to a Gaussian mixture, meaning in $\mathbf{B}_k X$ and with covariance matrix $\boldsymbol{\Sigma}_k= (\boldsymbol{P}_k^t \boldsymbol{P}_k)^{-1}=\sigma I_q$, for the cluster $k$.
We run our procedures with the number of components varying in $\mathcal{K} = \{2,\ldots,5 \}$.

The model $1$ illustrates our procedures in low dimension models.
Moreover, it is chosen in the next section to illustrate each step of the procedure (variable selection, models construction and model selection). 
Model $5$ is considered to be high-dimensional, because $p \times K > n$. 
The model $2$ is easier than the others, because clusters are not so closed to each other according to the noise. 
Model $3$ is constructed as the models $1$ and $2$, but $n$ is small and the noise is more important. We will see that the clustering is more difficult in that case.
Model $4$ has a larger SNR, nevertheless, the problem of clustering is difficult, because each $\mathbf{B}_k$ is closer to the others.

Our procedures are run $20$ times and we compute statistics on our results over those $20$ simulations: it is a small number, but the whole procedure is time-consuming and results are convincing enough.

For the initialization step, we repeat $50$ times the initialization and keep the one which maximizes the log-likelihood function after $10$ iterations.
Those choices are size-dependent, a numerical study which is not reported here concludes that it is enough in that case.
\subsection{Sparsity and model selection}
To illustrate the both procedures, all the analyses made in this section are done from the model $1$, since the choice of each step is clear.

\begin{figure}[!ht]
\centering
\begin{minipage}[c]{.46\linewidth}
       \includegraphics[scale=0.5]{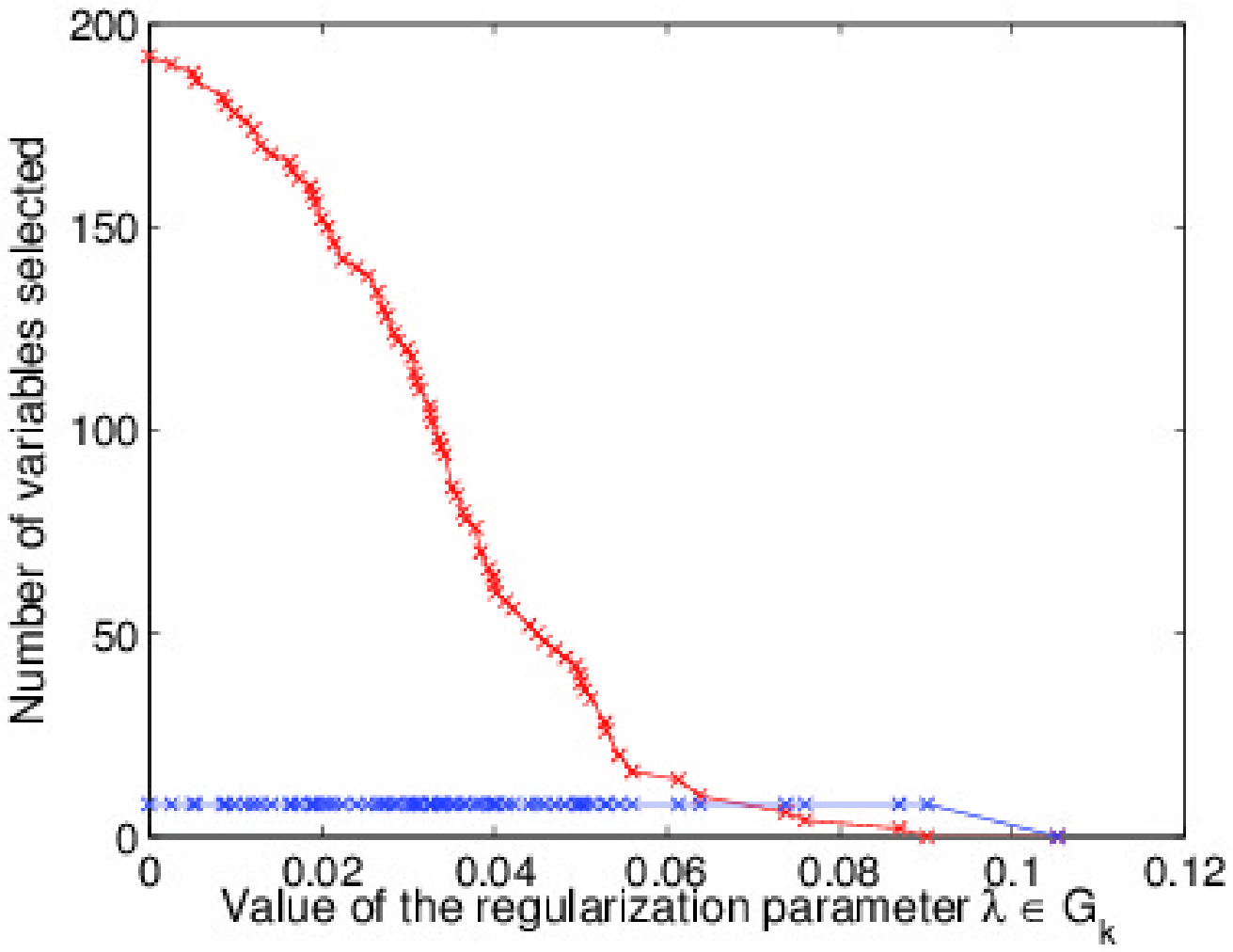}
 \caption[Number of FR and TR]{For one simulation, the number of false relevant (in red color) and true relevant (in blue color) variables generated by the Lasso, by varying the regularization parameter $\lambda$ in the grid $G_2$.}
 \label{selecVar}
   \end{minipage} \hfill
   \begin{minipage}[c]{.46\linewidth}
       \includegraphics[scale=0.5]{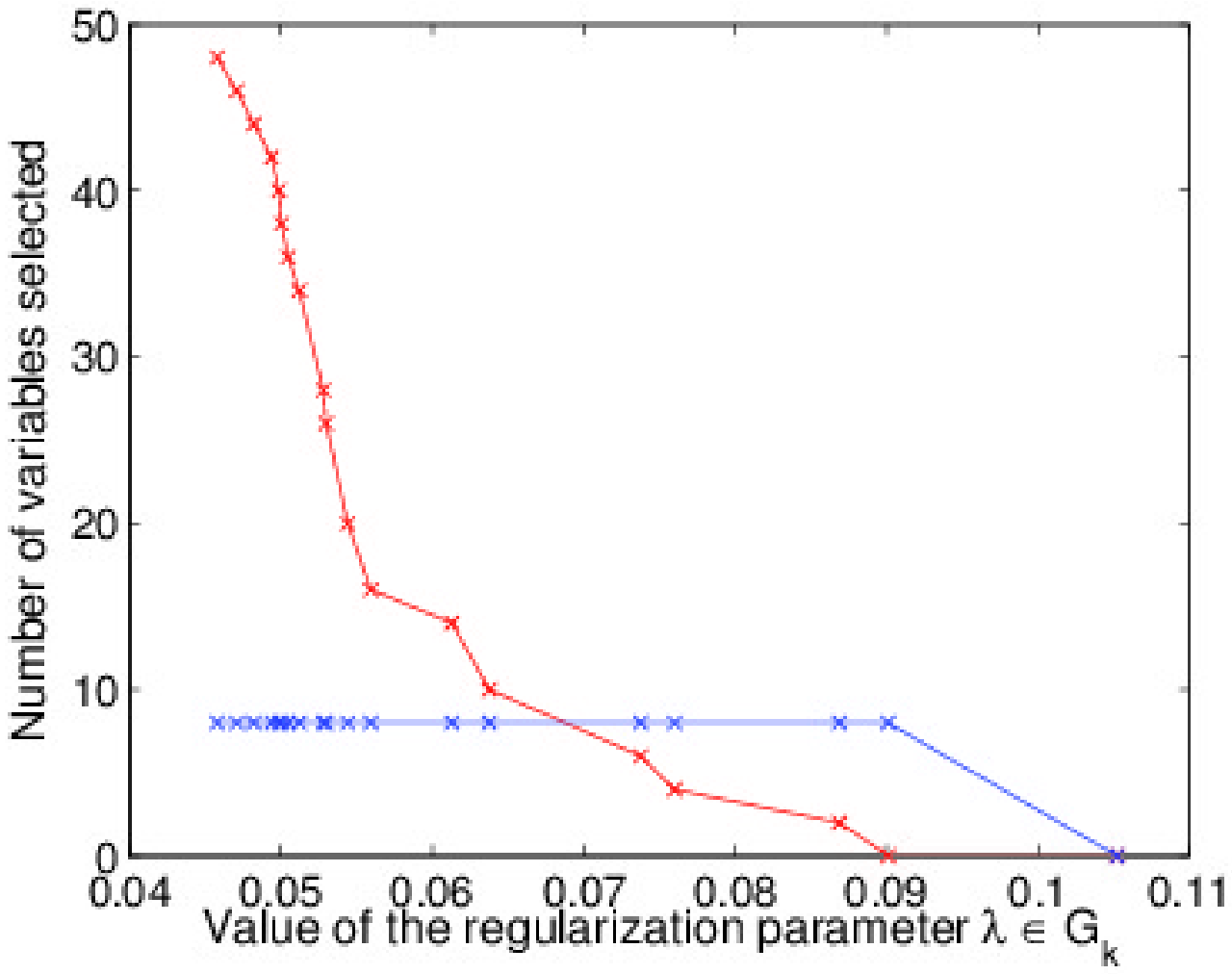}
       
 \caption[Zoom in on number of FR and Tr]{For one simulation, zoom in on the number of false relevant (in red color) and true relevant (in blue color) variables generated by the Lasso, by varying the regularization parameter $\lambda$ around the interesting values.}
 \label{selecVar2}
   \end{minipage}
 \end{figure}

Firstly, we compute the grid of regularization parameters.
More precisely, each regularization parameter is computed from maximum likelihood estimations (using EM algorithm) and gives an associated sparsity (computed by the Lasso estimator, using again the EM algorithm). 
In Figure \eqref{selecVar} and Figure \eqref{selecVar2}, the collection of relevant variables selected by this grid are plotted.

Firstly, we could notice that the number of relevant variables selected by the Lasso decreases with the regularization parameter. 
We could analyze more precisely which variables are selected, that is to say if we select true relevant or false relevant variables.
If the regularization parameter is not too large, the true relevant variables are selected. Even more, if the regularization parameter is well-chosen, 
we select only the true relevant variables.
In our example, we remark that if $\lambda=0.09$, we have selected exactly the true relevant variables.
This grid construction seems to be well-chosen according to these simulations.

From this variable selection, each procedure (Lasso-MLE or Lasso-Rank) leads to a model collection, varying the sparsity thanks to the regularization parameters grid and the number of components.

Among this collection, we select a model with the slope heuristic.

\begin{figure}
   \begin{minipage}[c]{.46\linewidth}
      \includegraphics[scale=0.5]{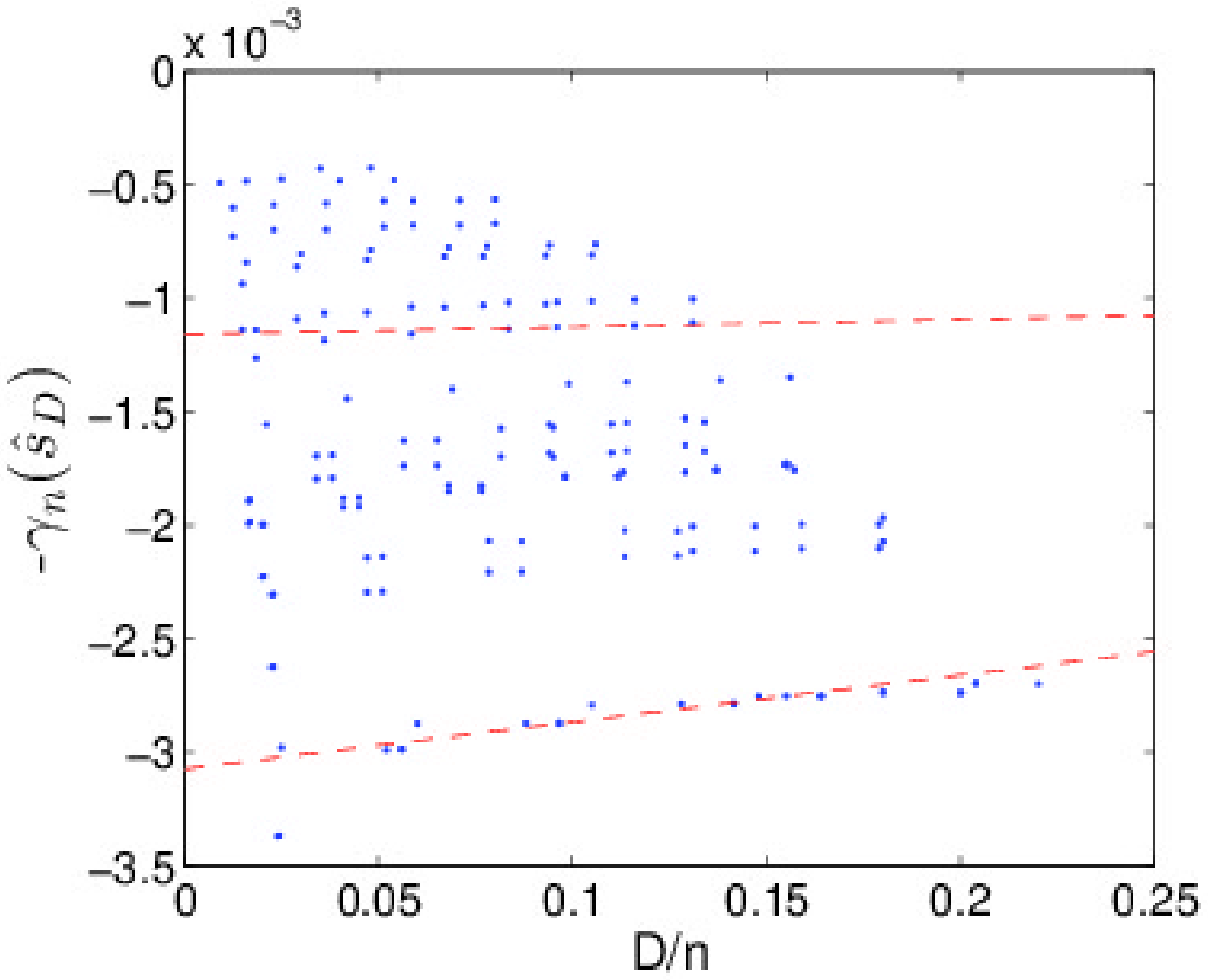}
      \caption[Slope graph for our Lasso-Rank procedure]{For one simulation, slope graph get by our Lasso-Rank procedure. For large dimensions, we observe a linear part.}
      \label{slope1}
   \end{minipage} \hfill
   \begin{minipage}[c]{.5\linewidth}
      \includegraphics[scale=0.5]{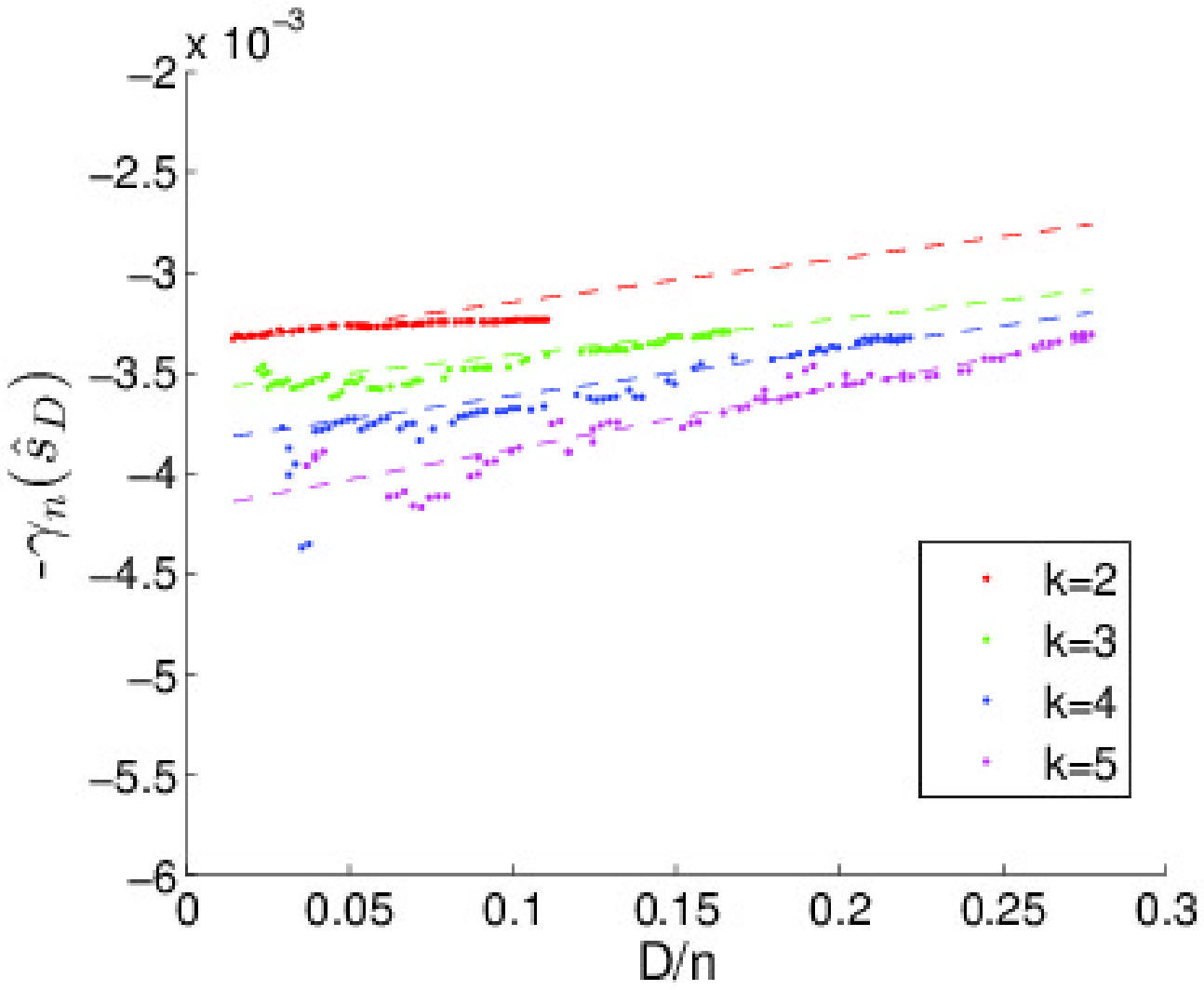}
      \caption[Slope graph for our Lasso-MLE procedure]{For one simulation, slope graph get by our Lasso-MLE procedure. For large dimensions, we observe a linear part.}
      \label{slope2}
   \end{minipage}
\end{figure}

We want to select the best model by improving a penalized criterion. This penalty is computed by performing a linear regression on the couples of points 
$\{ (D/n ; -1/n\sum_{i=1}^n \log (\hat{h}_{D}(\boldsymbol{y}_i|\boldsymbol{x}_i))) \}$. 
The slope $\hat{\kappa}$ allows us to have access to the best model, the one with dimension $\hat{D}$ minimizing 
$$-\frac{1}{n} \sum_{i=1}^n \log (\hat{h}_{D}(\boldsymbol{y}_i|\boldsymbol{x}_i))+2 \hat{\kappa} \frac{D}{n}.$$
In practice, we have to look if couples of points have a linear comportment.
For each procedure, we construct a different model collection and we have to justify this behavior.
Figures \eqref{slope1} and \eqref{slope2} represent the log-likelihood in function of the dimension of the models, for model collections constructed respectively by the Lasso-MLE procedure and by the Lasso-Rank procedure.
The couples are plotted by points, whereas the estimated slope is specified by a dotted line.
We could observe more than a line ($4$ for the Lasso-MLE procedure, more for the Lasso-Rank procedure).
This phenomenon could be explained by a linear behavior for each mixture, fixing the number of classes and for each rank.
Nevertheless, slopes are almost the same and select the same model. In practice, we estimate the slope with the Capushe package described in \cite{baudry}.

\subsection{Assessment}
We compare our procedures to three other procedures on the simulated models $1$, $2$, $3$ and $4$. 

Firstly, let us give some remarks about the model $1$.
For each procedure, we get a good clustering and a very low Kullback-Leibler divergence. 
Indeed, the sample size is large and the estimations are good. That is the reason why, in this subsection, we focus on models $2$, $3$ and $4$.

To compare our procedures with others, the Kullback-Leibler  divergence with the true density and the ARI (the Adjusted Rand Index, measuring the similarity between two data clusterings, knowing that the closer to $1$ the ARI, the more similar the two partitions) are computed, 
and we note which variables are selected and how many clusters are selected. For more details on the ARI, see \cite{Arabie}.

From the Lasso-MLE model collection, we construct two models, to compare our procedures with. 
We compute the oracle (the model which minimizes the Kullback-Leibler divergence with the true density) and the model which is selected by the BIC criterion instead of the slope heuristic.
According to the oracle model, we know how good we could get from this model collection for the Kullback-Leibler divergence and how this model, 
as good it is possible for the log-likelihood, performs the clustering.

The third procedure we compare with is the maximum likelihood estimator, assuming that we know how many clusters there are, fixed to $2$.
We use this procedure to show that variable selection is needed.

In each case, we apply the MAP principle, to compare clusterings.

We do not plot the Kullback-Leibler divergence for the MLE procedure, because values are too high and make the boxplots unreadable.

\begin{figure}
   \begin{minipage}[c]{.46\linewidth}
   \centering
      \includegraphics[scale=0.45]{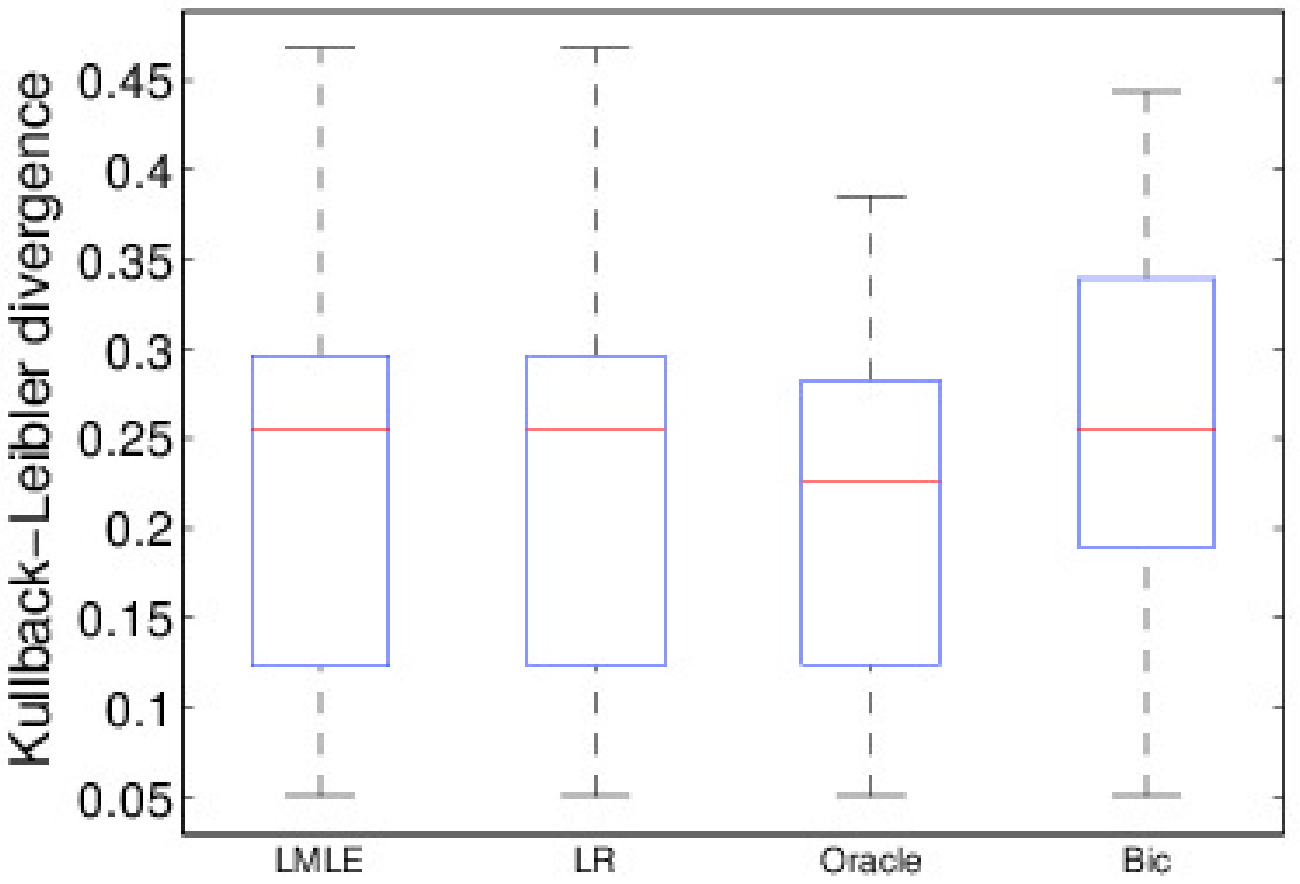}
\caption[Boxplot of the Kullback-Leibler divergence] {Boxplots of the Kullback-Leibler divergence between the true model and the one selected by the procedure over the $20$ simulations, for the Lasso-MLE procedure (LMLE), the Lasso-Rank procedure (LR), the oracle (Oracle), the BIC estimator (BIC) for the model $2$.}
\label{KLsec1}
   \end{minipage} \hfill
   \begin{minipage}[c]{.46\linewidth}
   \centering
   \includegraphics[scale=0.45]{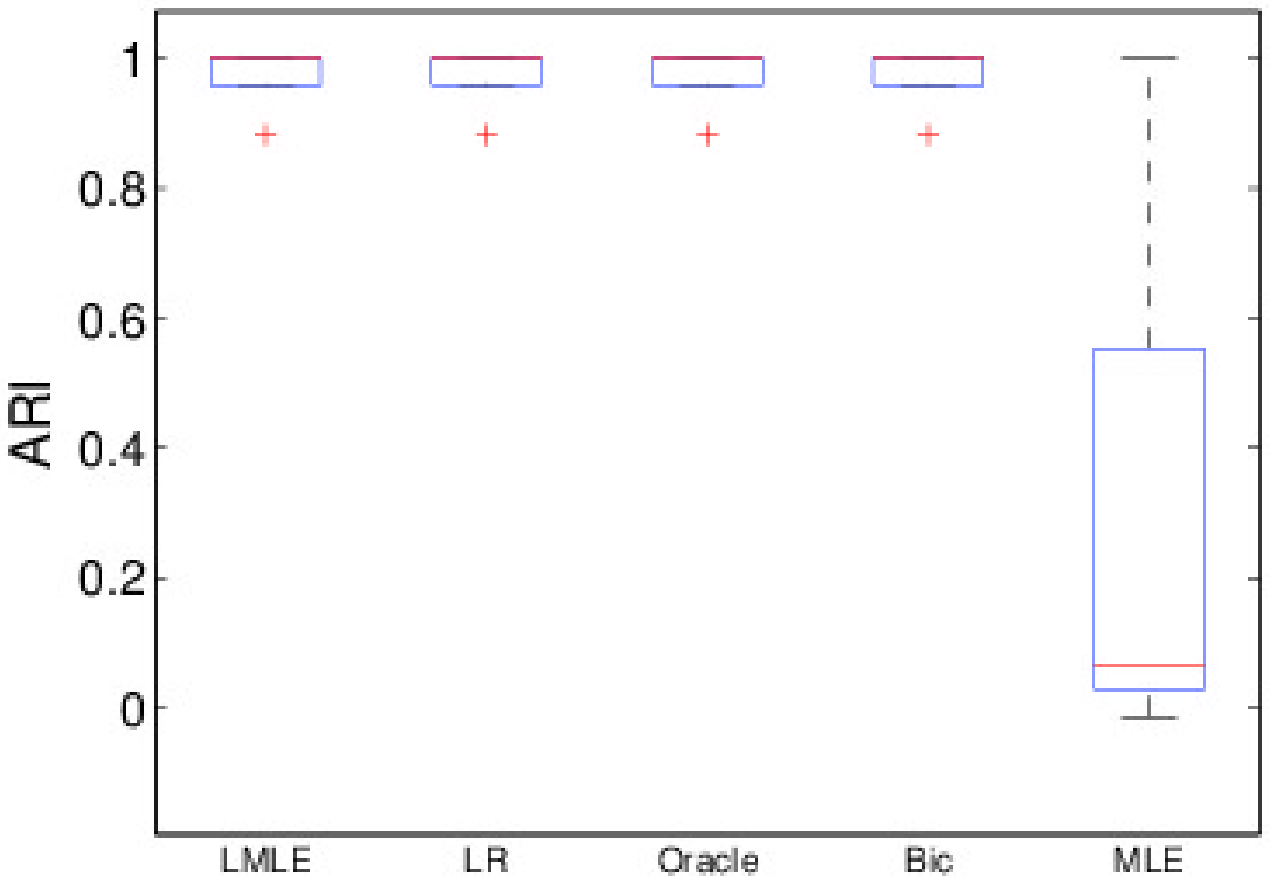}
\caption[Boxplot of the ARI]{Boxplots of the ARI over the $20$ simulations, for the Lasso-MLE procedure (LMLE), the Lasso-Rank procedure (LR), the oracle (Oracle), the BIC estimator (BIC) and the MLE (MLE) for the model $2$.}
\label{ARIsec1}
   \end{minipage}
\end{figure}

\begin{figure}
   \begin{minipage}[c]{.46\linewidth}
\centering
   \includegraphics[scale=0.45]{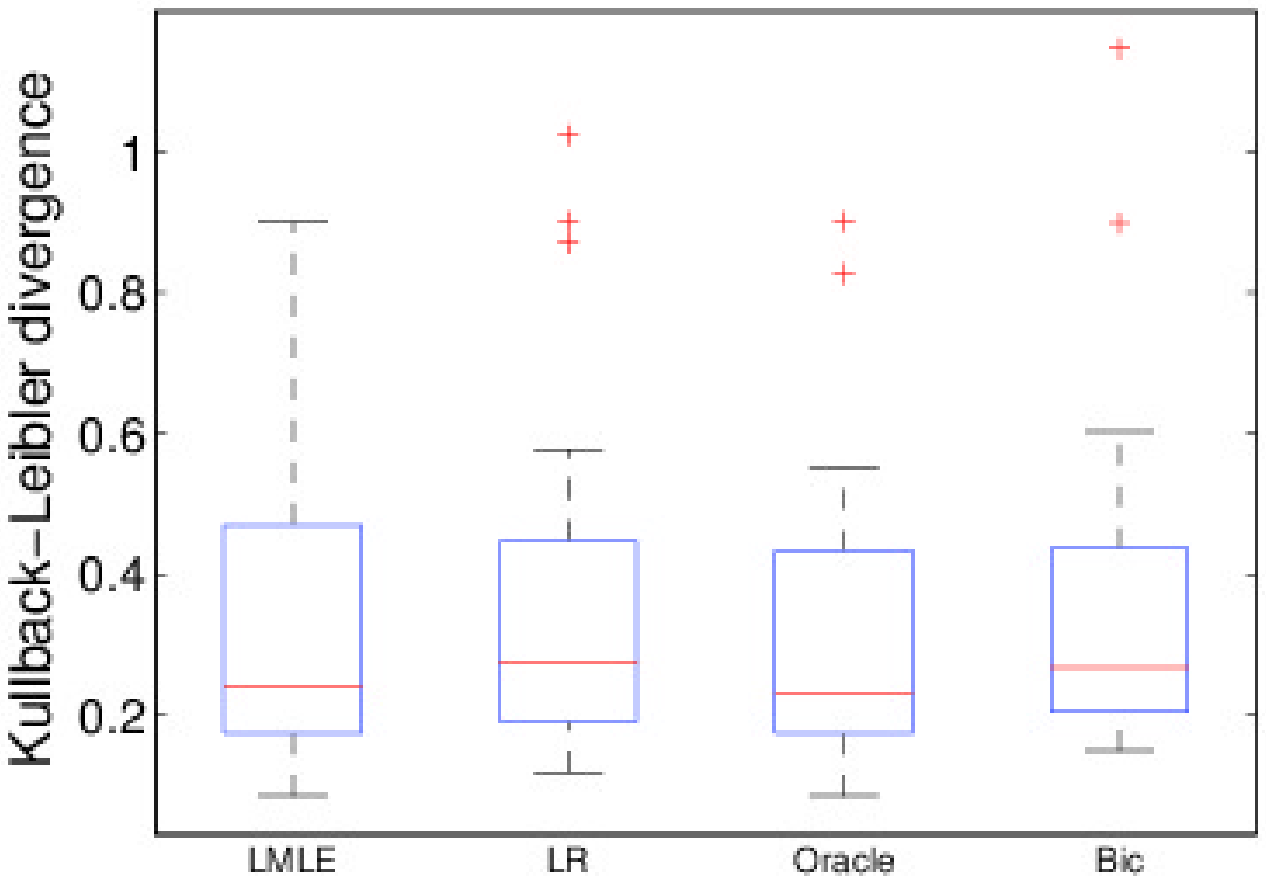}
\caption[Boxplot of the Kullback-Leibler divergence]{Boxplots of the Kullback-Leibler divergence between the true model and the one selected by the procedure over the $20$ simulations, for the Lasso-MLE procedure (LMLE), the Lasso-Rank procedure (LR), the oracle (Oracle), the BIC estimator (BIC) for the model $3$.}
\label{KL2sec1}
   \end{minipage} \hfill
   \begin{minipage}[c]{.46\linewidth}
\centering
   \includegraphics[scale=0.45]{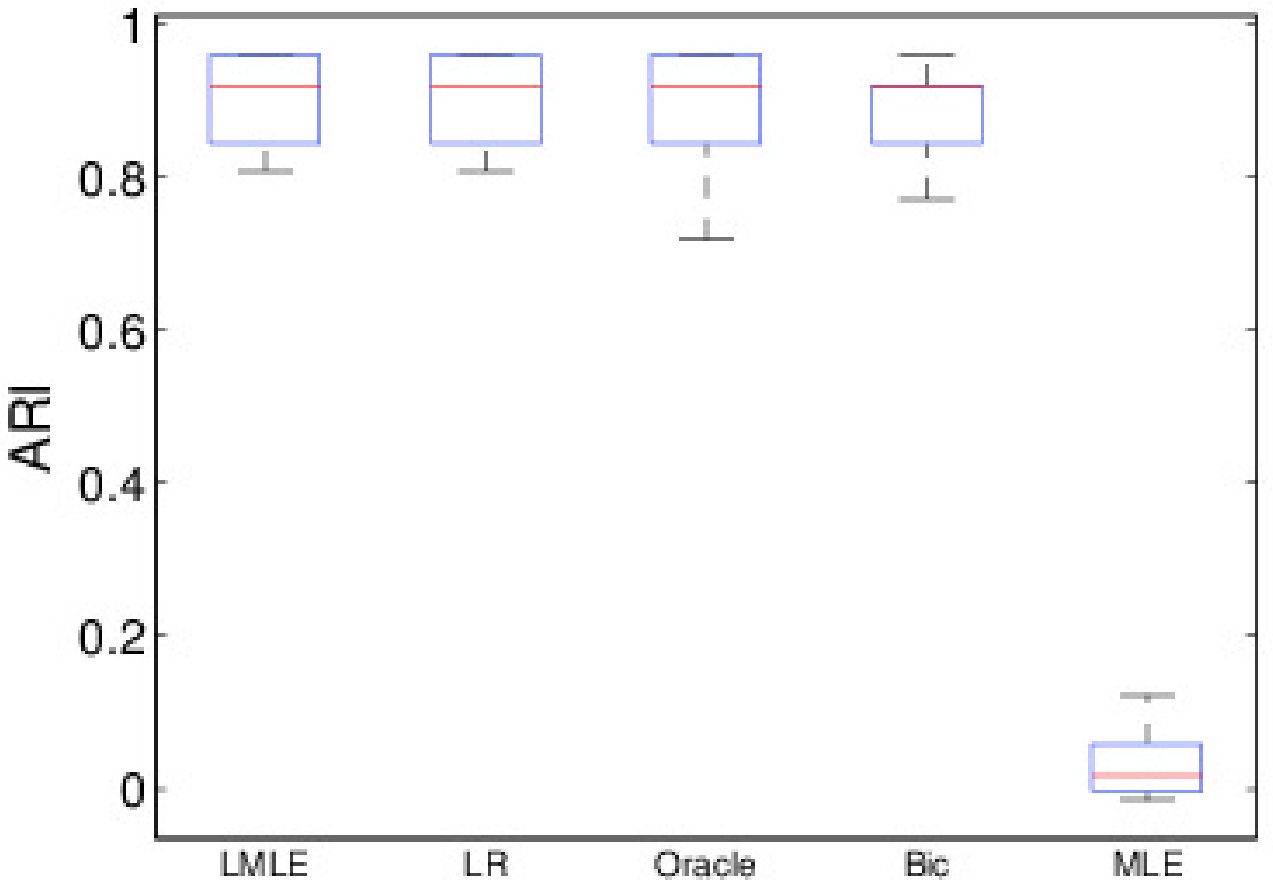}
\caption[Boxplot of the ARI]{Boxplots of the ARI over the $20$ simulations, for the Lasso-MLE procedure (LMLE), the Lasso-Rank procedure (LR), the oracle (Oracle), the BIC estimator (BIC) and the MLE (MLE) for the model $3$.}
\label{ARI2sec1}
   \end{minipage}
\end{figure}
\begin{figure}
   \begin{minipage}[c]{.46\linewidth}
\centering
   \includegraphics[scale=0.45]{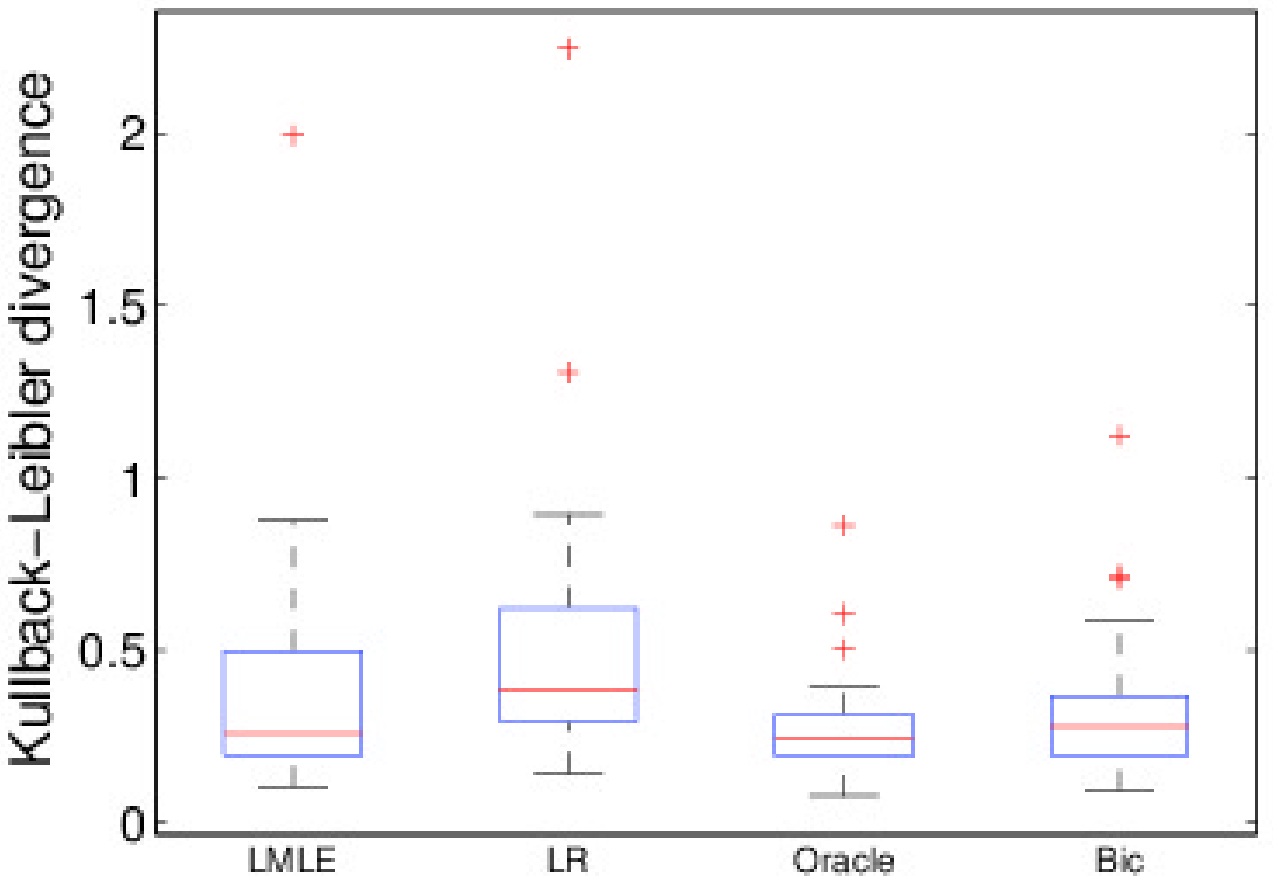}
\caption[Boxplot of the Kullback-Leibler divergence]{Boxplots of the Kullback-Leibler divergence between the true model and the one selected by the procedure over the $20$ simulations, for the Lasso-MLE procedure (LMLE), the Lasso-Rank procedure (LR), the oracle (Oracle), the BIC estimator (BIC) for the model $4$.}
\label{KL3sec1}
   \end{minipage} \hfill
   \begin{minipage}[c]{.46\linewidth}
\centering
   \includegraphics[scale=0.4]{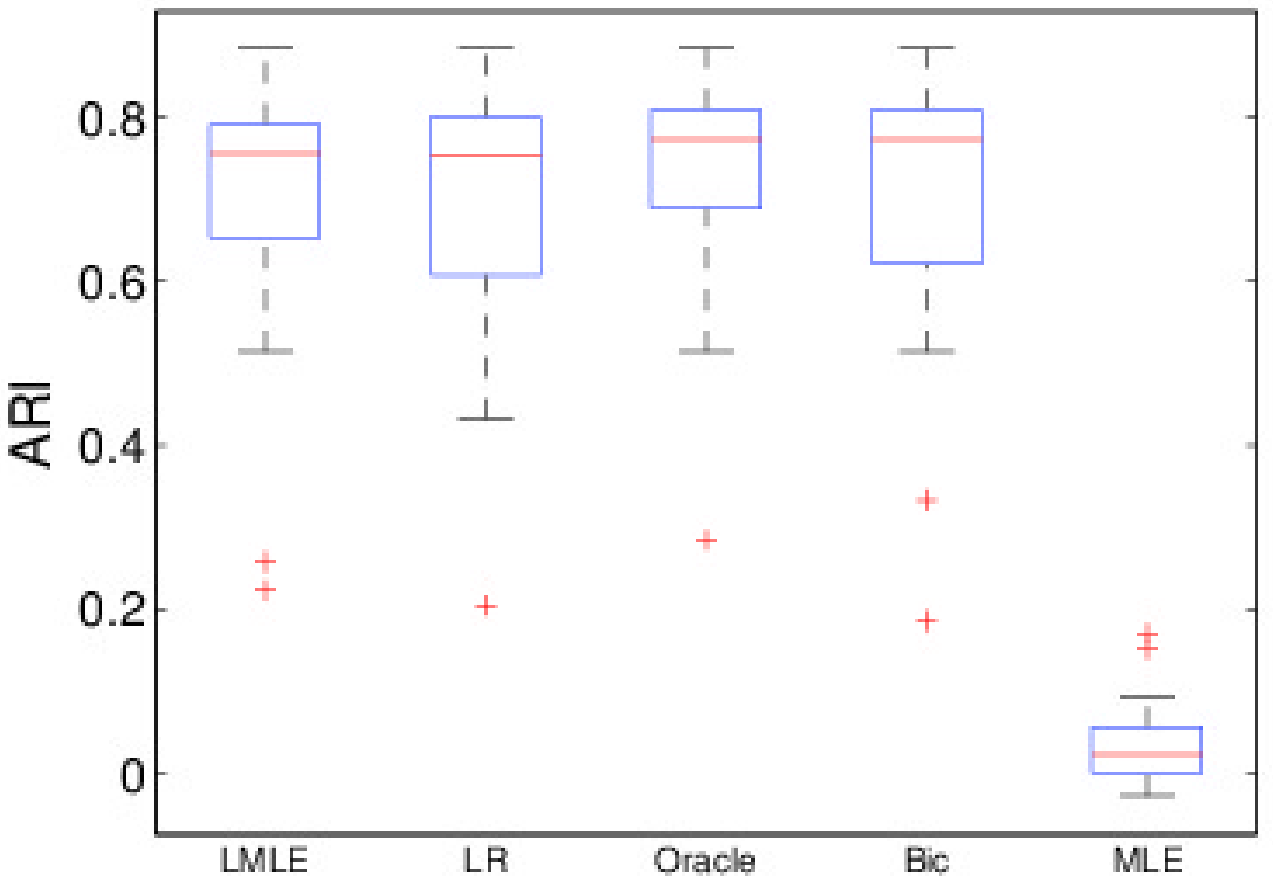}
\caption[Boxplot of the ARI]{Boxplots of the ARI over the $20$ simulations, for the Lasso-MLE procedure (LMLE), the Lasso-Rank procedure (LR), the oracle (Oracle), the BIC estimator (BIC) and the MLE (MLE) for the model $4$.}
\label{ARI3sec1}
   \end{minipage}
\end{figure}
\begin{figure}
   \begin{minipage}[c]{.46\linewidth}
\centering
   \includegraphics[scale=0.4]{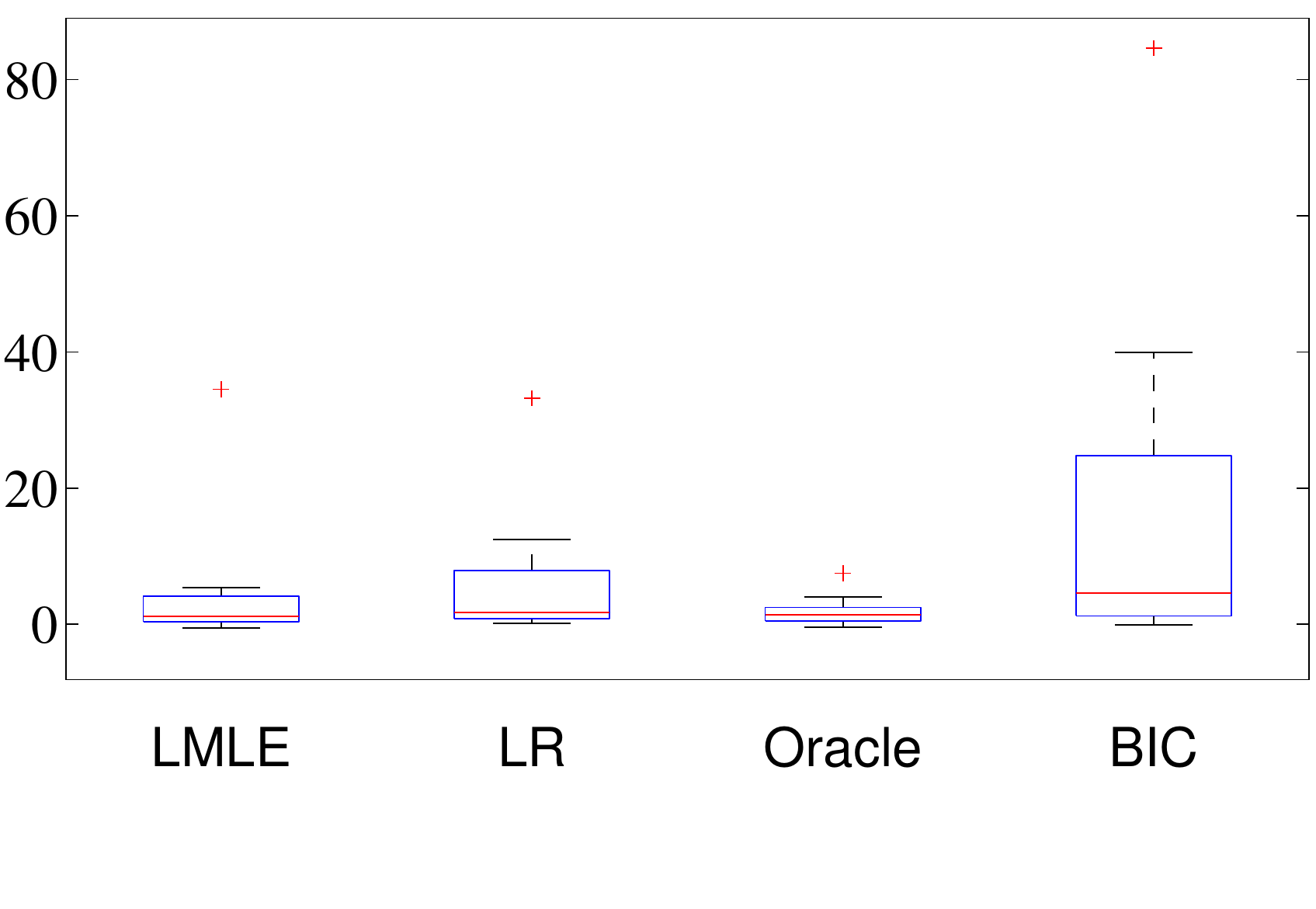}
\caption[Boxplot of the Kullback-Leibler divergence]{Boxplots of the Kullback-Leibler divergence between the true model and the one selected by the procedure over the $20$ simulations, for the Lasso-MLE procedure (LMLE), the Lasso-Rank procedure (LR), the oracle (Oracle), the BIC estimator (BIC) for the model $5$.}
\label{KL5sec1}
   \end{minipage} \hfill
   \begin{minipage}[c]{.46\linewidth}
\centering
   \includegraphics[scale=0.33]{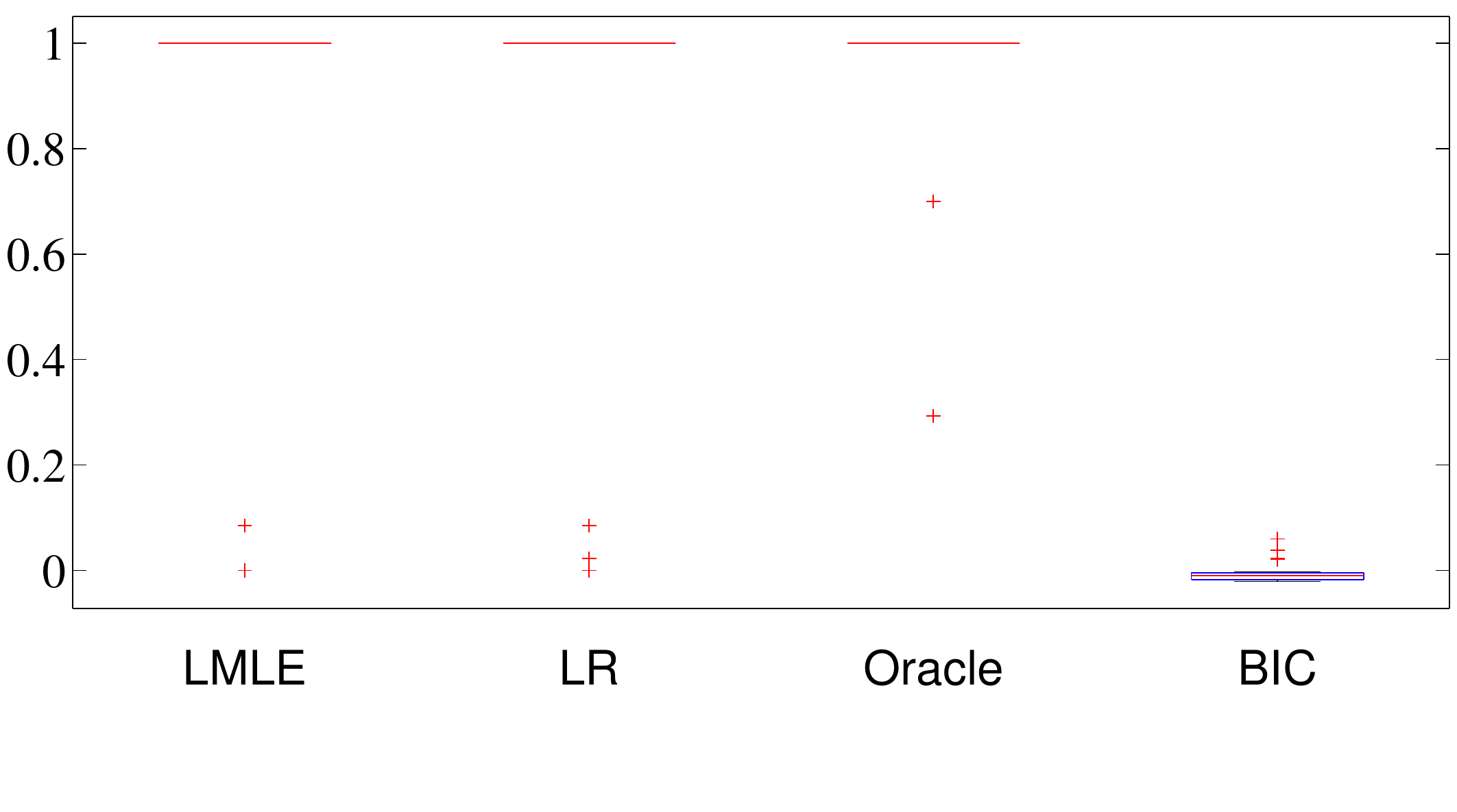}
\caption[Boxplot of the ARI]{Boxplots of the ARI over the $20$ simulations, for the Lasso-MLE procedure (LMLE), the Lasso-Rank procedure (LR), the oracle (Oracle), the BIC estimator (BIC) and the MLE (MLE) for the model $5$.}
\label{ARI5sec1}
   \end{minipage}
\end{figure}
For the model $2$, according to Figure \eqref{KLsec1} for the Kullback-Leibler divergence and Figure \eqref{ARIsec1} for the ARI, the Kullback-Leibler divergence is small and the ARI is close to $1$, except for the MLE procedure.
Boxplots are still readable with those values, but it is important to highlight that variable selection is needed, even in a model with reasonable dimension.
The model collections are then well constructed.
The model $3$ is more difficult, because the noise is higher. That is why results, summarized in Figures \eqref{KL2sec1} and \eqref{ARI2sec1}, are not as good as for the model $2$. Nevertheless, our procedures lead to the best ARI and the Kullback-Leibler divergences are close to the one of the oracle.
We could make the same remarks for the model $4$. In this study, means are closer, according to the noise. Results are summarized in Figures \eqref{KL3sec1} and \eqref{ARI3sec1}.
The model $5$ is in high-dimension. Models selected by the BIC criterion are bad, in comparison with models selected by our procedures, or with oracles.
They are bad for estimation, according to the Kullback-Leibler divergence boxplots in Figure \eqref{KL5sec1}, but also for clustering, according to Figure \eqref{ARI5sec1}.
Our models are not as well as constructed as previously, but it is explained by the high-dimensional context.
It is explained by the high Kullback-Leibler divergence. However, our performances in clustering are really good.

Note that
the Kullback-Leibler divergence is smaller for the Lasso-MLE procedure, thanks to the maximum likelihood refitting.
Moreover, the true model has not any matrix structure.
If we compare with the MLE, where we do not use the sparsity hypothesis, we could conclude that estimations are not satisfactory, which could be
explained by an high-dimensional issue.
The Lasso-MLE procedure, the Lasso-Rank procedure and the BIC model work almost as well as the oracle, which mind that the models are well selected.

\begin{table}[H]
\centering
\begin{tabular}{|c||cc||cc||cc|}
 \hline
 && Model 2&&Model 3&&Model 4\\
 \hline
 \hline
 Procedure &  TR&FR  &TR&FR &TR&FR\\ \hline 
  Lasso-MLE & $ 8 (0)$&$ 2.2 (6.9) $& $8(0) $& $4.3 (28.8)$ &$8(0)$&$2(13,5)$\\ \hline
  Lasso-Rank&  $8 (0)$&$24 (0)$ & $8 (0)$& $24(0)$& $8(0)$&$24(0)$\\ \hline
  Oracle&  $8 (0)$&$ 1.5 (3.3)$& $7.8 (0.2)$&$2.2 (11.7)$ & $8 (0)$&$ 0.8 (2.6)$\\ \hline
  BIC estimator& $8 (0)$&$ 2.6 (15.8)$ & $7.8 (0.2)$& $5.7 (64.8)$ &$ 8 (0)$&$2.6 (11.8)$\\ \hline
\end{tabular}
\caption[Mean number of TR and FR]{Mean number \{TR, FR\} of true relevant and false relevant variables over the $20$ simulations for each procedure, for models $2$, $3$ and $4$. The standard deviations are put into parenthesis.}
\label{varSelect}
\end{table}

In Table \eqref{varSelect}, we summarize results about variable selection. For each model and for each procedure, we compute how many true relevant and false relevant variables are selected.

The true model has $8$ relevant variables, which are always recognized.
The Lasso-MLE has less false relevant variables than the others, which means that the true structure was found.
The Lasso-Rank has $24$ false relevant variables, because of the matrix structure: the true rank in each component was $4$, then the estimator restricted on relevant variables is a $4 \times 4$ matrix and we get $12$ false relevant variables in each component.
Nevertheless, we do not have more variables, that is to say the model constructed is the best possible.
The BIC estimator and the oracle have a large variability for the false relevant variables.

For the number of components, we find that all the procedures have selected the true number $2$.
 
 Thanks to the MLE, the first procedure has good estimations (better than the second one). However, depending on the data, the second procedure could be more attractive.
 If there is a matrix structure, for example if most of the variation of the response $\boldsymbol{Y}$ is caught by a small number of linear combinations of the predictors, the second procedure will work better.
 
We could conclude that the model collection is well constructed and that the clustering is well-done.

  \section{Functional dataset}
  \label{FunctionalDatasets}
  One of the main interest of our methods is to be applied to functional dataset.
  Indeed, in different fields of applications, the considered data are functions.
The functional data analysis has been popularized first by \cite{Ramsay}. It gives a description of the main tools to analyze functional dataset.
Another book is the one of \cite{Ferraty}.
In a density estimation framework, we can refer to \cite{GarethSugar}, who used a model-based approach for clustering functional data.
About functional regression, however, the main part of the existing literature is focused on $\boldsymbol{Y}$ scalar and $\boldsymbol{X}$ functional. For example, 
we can refer to \cite{ZhaoOgdenReiss} for using wavelet basis in linear model, \cite{Yao} for functional mixture regression,
or \cite{Ciarleglio} for using wavelet basis in functional mixture regression.
In this section, we focus on $\boldsymbol{Y}$ and $\boldsymbol{X}$ both functional.
In this regression context, we could be interested in clustering: it leads to identify individuals involved in the same relation between $\boldsymbol{Y}$ and $\boldsymbol{X}$.
Denote that, with functional dataset, we have to denoise and smooth signals to remove the noise and capture only the important patterns in the data.
Here, we explain how our procedures can be applied in this context. Note that we could apply our procedures with scalar response and functional 
regressor, or, on the contrary, with functional response and scalar regressor. We explain how our procedures
work in the more general case.
Remark that we focus here on the wavelet basis, to take advantage of the time-scale decomposition, but the same analysis is available
on Fourier basis or splines.
\label{funct}
\subsection{Functional regression model}
 Suppose we observe a centered sample of functions $(\boldsymbol{f}_i,\boldsymbol{g}_i)_{1 \leq i \leq n}$, associated with the random variables $(\boldsymbol{F},\boldsymbol{G})$, 
 coming from a probability distribution with unknown conditional density $s^{*}$. 
We want to estimate this model by a functional mixture model: if the variables $(\boldsymbol{F},\boldsymbol{G})$ come from the component $k$, there exists a function $\mathbf{B}_k$ such that
\begin{equation}
\label{modelefonc}
\boldsymbol{G}(t)= \int_{I_x}\boldsymbol{F}(u) \mathbf{B}_k(u,t) du + \boldsymbol{\varepsilon}(t),
\end{equation}
where $\boldsymbol{\varepsilon}$ is a residual function. This linear model is introduced in \cite{Ramsay}.  They propose to project onto basis the response and the regressors.
We extend their model in mixture model, to consider several subgroups for a sample.

If we assume that, for all $t$, for all  $i \in \{1,\ldots,n\}$, $\boldsymbol{\varepsilon}_i(t) \sim \mathcal{N}(0,\boldsymbol{\Sigma}_k)$, the model \eqref{modelefonc} is an integrated version of the model \eqref{modeleFMR}.
Depending on the cluster $k$, the linear relation of $\boldsymbol{G}$ with respect to $\boldsymbol{F}$ is described by the function $\mathbf{B}_k$.
\subsection{Two procedures to deal with functional dataset}
\subsubsection{Projection onto a wavelet basis}

 To deal with functional data, we project them onto some basis, to obtain data as described in the 
 Gaussian mixture regression models \eqref{modeleFMR}.
 In this article, we choose to deal with wavelet basis, given that they represent localized features of functions in a sparse way. 
 If the coefficient matrices ${\textbf{x}}$ and ${\textbf{y}}$ are sparse,
 the regression matrix $\mathbf{B}$ has more chance to be sparse.
 Moreover, we could represent a signal with few coefficients, which reduces the dimension.
 For details about the wavelet theory, see \cite{Mallat}.
 
Let us begin with an overview of some important aspects of wavelet basis.

Let $\psi$ be a real wavelet function, satisfying
$$ \psi \in L^1 \cap L^2, t \psi \in L^1, \text{ and } \int_{\mathbb{R}} \psi(t) dt =0.$$
We denote by $\psi_{l,h}$ the function defined from $\psi$ by dyadic dilation and translation:
$$\psi_{l,h}(t)=2^{l/2} \psi(2^l t-h) \text{ for } (l,h) \in \mathbb{Z}^2.$$
We could define the wavelet coefficients of a signal $f$ by
$$d_{l,h}(f) =\int_{\mathbb{R}} f(t) \psi_{l,h}(t)dt \text{ for } (l,h) \in \mathbb{Z}^2.$$
Let $\varphi$ be a scaling function related to $\psi$ and $\varphi_{l,h}$ be the dilatation and translation of $\varphi$ for $(l,h) \in \mathbb{Z}^2$. We also define, for $(l,h) \in \mathbb{Z}^2$,
$\mathbf{B}_{l,h}(f) = \int_{\mathbb{R}} f(t) \varphi_{l,h}(t)dt$.

Note that scaling functions serve to construct approximations of the function of interest, while the wavelet functions serve to provide the details not captured by successive approximations.

We denote by $V_l$ the space generated by $\{\varphi_{l,h}\}_{h \in \mathbb{Z}}$ and by $W_l$ the space generated by $\{ \psi_{l,h}\}_{h \in \mathbb{Z}}$ for all $l \in \mathbb{Z}$.
Remark that
\begin{align*}
 V_{l-1} &= V_l \oplus W_l \text{ for all } l \in \mathbb{Z} ;\\
 L^2 &= \oplus_{l \in \mathbb{Z}} W_l.
\end{align*}

Let $L \in \mathbb{N}^*$. For a signal $f$, we could define the approximation at the level $L$ by
$$A_L = \sum_{l >L} \sum_{h \in \mathbb{Z}} d_{l,h} \psi_{l,h} ;$$
and $f$ could be decomposed by the approximation at the level $L$ and the details $(d_{l,h})_{l <L}$.

The decomposition of the basis between the scaling function and the wavelet function emphasizes the local nature of the wavelets.
It is an important aspect in our procedures, because we want to know which details allow us to cluster two observations together.

We consider the sample $(\boldsymbol{f}_i,\boldsymbol{g}_i)_{1\leq i \leq n}$ and introduce the wavelet expansion of $\boldsymbol{f}_i$ in the basis $\mathcal{B}$: for all $t \in [0,1] $,
$$\boldsymbol{f}_i(t)= \underbrace{\sum_{h \in \mathbb{Z}} \mathbf{B}_{L,h} (\boldsymbol{f}_i) \varphi_{L,h}(t)}_{A_L} + \sum_{l\leq L} \sum_{h\in \mathbb{Z}}d_{l,h}(\boldsymbol{f}_i) \psi_{l,h}(t).$$
The collection $\{\mathbf{B}_{L,h}(\boldsymbol{f}_i),d_{l,h}(\boldsymbol{f}_i) \}_{l \leq L,h \in \mathbb{Z}}$ is the Discrete Wavelet Transform (DWT) of $f$ in the basis $\mathcal{B}$.

Because we project onto an orthonormal basis, this leads to a $n$-sample $(\boldsymbol{x}_1,\ldots,\boldsymbol{x}_n)$ of wavelet coefficient decomposition vectors, with
$$\boldsymbol{f}_i=\boldsymbol W \boldsymbol{x}_i ;$$
in which $\boldsymbol{x}_i$ is the vector of the discretized values of the signal, $\boldsymbol{x}_i$ the matrix of coefficients in the basis $\mathcal{B}$, 
and $\boldsymbol W$ a $p \times p$ matrix defined by $\varphi$ and $\psi$.
The DWT can be performed by a computationally fast pyramid algorithm (see  \cite{Mallat2}).
In the same way, there exists $\boldsymbol W^{'}$ such that $\boldsymbol{g}_i=\boldsymbol W^{'} \boldsymbol{y}_i$, with ${\mathbf{y}}=(\boldsymbol{y}_1,\ldots,\boldsymbol{y}_n)$ a $n$ sample of wavelet coefficient decomposition vectors.
As the matrices $\boldsymbol W$ and $\boldsymbol W^{'}$ are orthogonal, we keep the mixture structure and the noise is also Gaussian.
We could consider the wavelet coefficient dataset $(\mathbf{{x}},\mathbf{{y}})=((\boldsymbol{x}_1,\boldsymbol{y}_1),\ldots,(\boldsymbol{x}_n,\boldsymbol{y}_n))$, which defines $n$ observations whose probability distribution could be modeled by the finite Gaussian mixture regression model \eqref{modeleFMR}.

\subsubsection{Our procedures}
We could apply our both procedures to this dataset and obtain a clustering of the data.
Indeed, rather than considering $(\mathbf{f},\mathbf{g})$, we run our procedures on the sample $(\mathbf{x},\mathbf{y})$, varying the number of clusters in $\mathcal{K}$.

The notion of relevant variable is natural: the function $\varphi_{l,h}$ or $\psi_{l,h}$ is irrelevant if it appears in none of the wavelet coefficient decomposition of the functions in each cluster.

\subsection{Numerical experiments}
We will illustrate our procedures on functional dataset by using the Matlab wavelet toolbox (see  \cite{poggi} for details). 
Firstly, we simulate functional dataset, where the true model belongs to the model collection.
Then, we run our procedure on an electricity dataset, to cluster successive days. We have access to time series, measured every half-hour, of a load consumption, on $70$ days.
We extract the signal of each day and construct couples by each day and its eve and we aim at clustering these couples.
To finish, we test our procedures on the well-known Tecator dataset. 
For each meat sample the data consists of a 100 channel spectrum of 
absorbances and the contents of fat. 
These experiments illustrate different aspects of our procedures.
Indeed, the simulated example proves that our procedures work in a functional context. The second example is a toy example used to validate the classification, on real data already studied and in which we clearly understand the clusters.
The last example illustrates the use of the classification to perform prediction and the description given by our procedures to the model constructed.
\subsubsection{Simulated functional data}
Firstly, we simulate a mixture regression model.
Let $\mathbf{f}$ be a sample of the noised cosine function, discretized on a $15$ points grid.
Let $\mathbf{g}$ be, depending on the cluster, either $\mathbf{f}$, or the function $-\mathbf{f}$, computed by a white-noise.

We use the Daubechies-2 basis at level $2$ to decompose the signal.

Our procedures are run $20$ times and the number of clusters are fixed to $\mathcal{K}=2$.
We construct a model collection by varying the sparsity for the Lasso-MLE procedure and by varying the sparsity and the ranks vector for the conditional mean for the Lasso-Rank procedure.
Among this model collection, we select a model with the slope heuristic.
We compare the models get by our procedures with the oracle model among the collection constructed by the Lasso-MLE procedure and with the model 
selected by the BIC criterion among this collection.
\begin{figure}[H]
      \centering
      \includegraphics[scale=0.5]{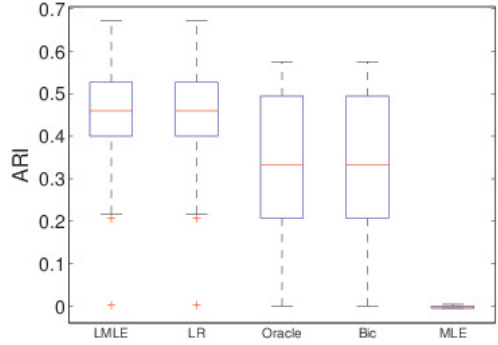}
      \caption[Boxplot of the ARI]{Boxplots of the ARI over the $20$ simulations, for the Lasso-MLE procedure (LMLE), the Lasso-Rank procedure (LR), the oracle (Oracle), the BIC estimator (Bic) and the MLE (MLE).}
      \label{ARIfonct}
\end{figure}
The ARIs are lower for models constructed by our models collection.
This simulated dataset proves that our procedures also perform clustering functional data, considering the projection dataset.
\subsubsection{Electricity dataset}
\label{Electricity}
\begin{figure}[H]
   \begin{minipage}[c]{.46\linewidth}
   \centering
\includegraphics[scale=0.3]{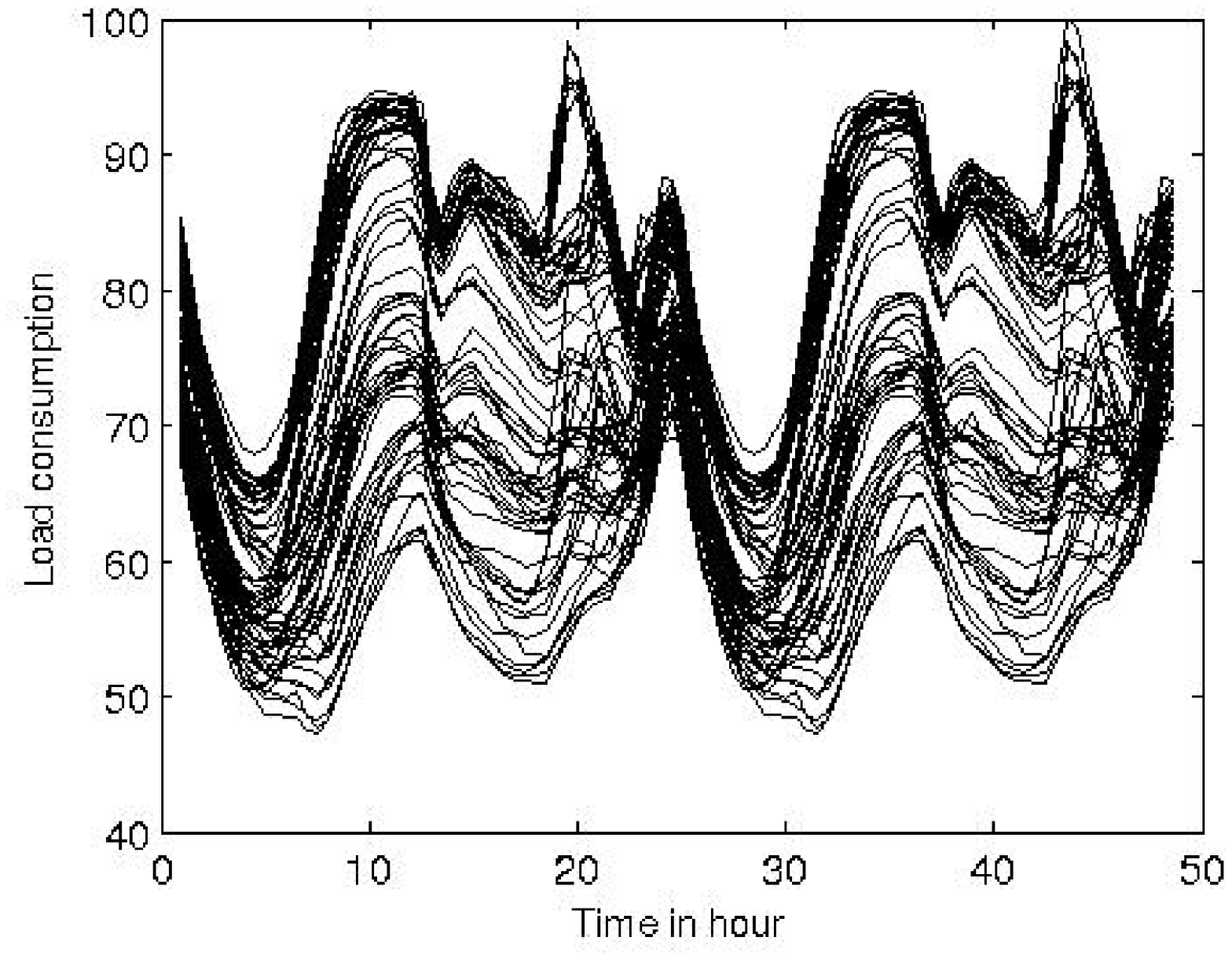}
  \caption{Plot of the $70$-sample of half-hour load consumption, on the two days.}
  \label{sample}
   \end{minipage} \hfill
   \begin{minipage}[c]{.46\linewidth}
   \centering
  \includegraphics[scale=0.2]{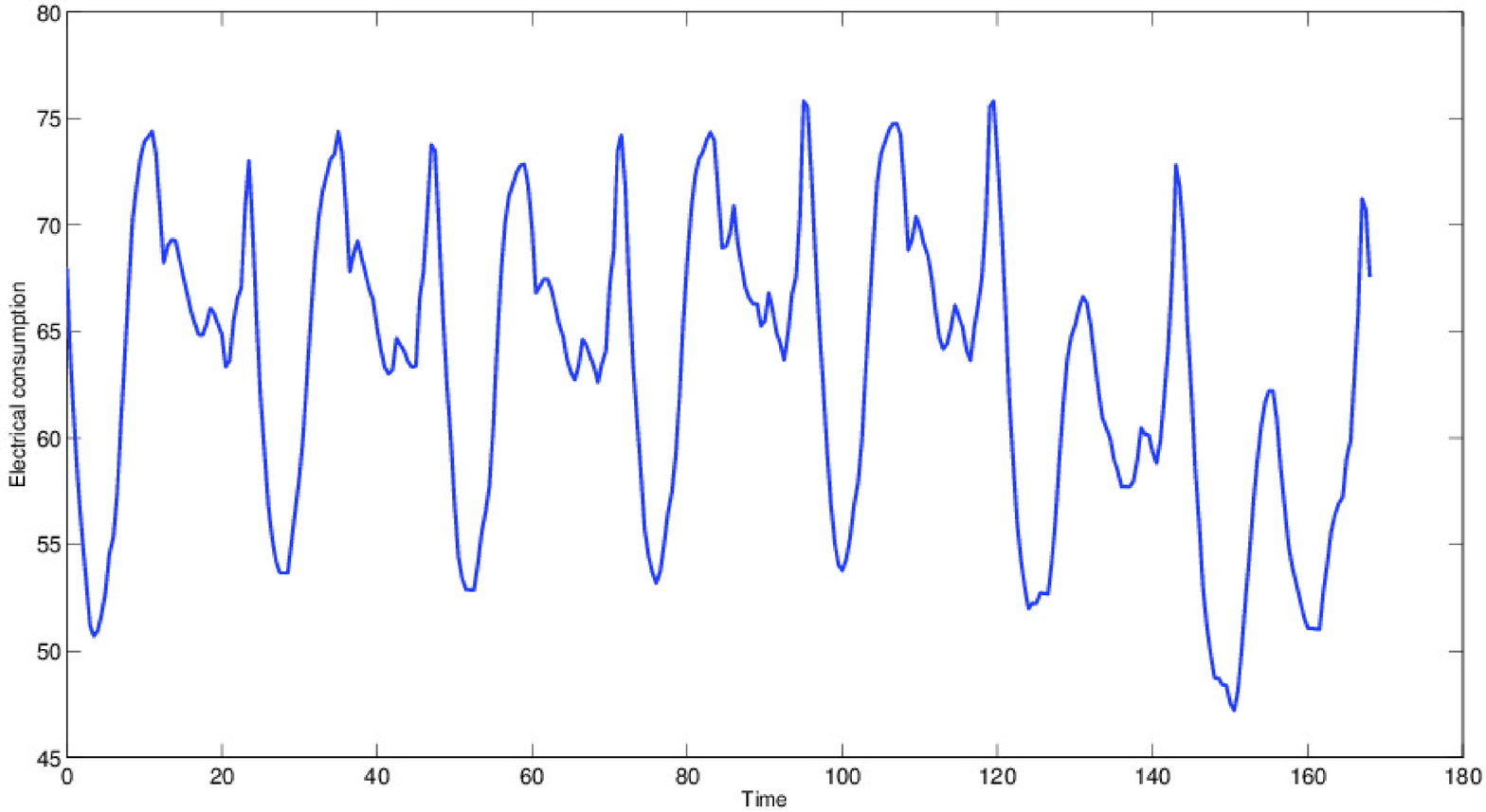}
  \caption{Plot of a week of load consumption.}
   \label{uneSemaine}
   \end{minipage}
\end{figure}

We also study the clustering on electricity dataset.
This example is studied in \cite{Misiti}.
We work on a sample of size $70$ of couples of days, which is plotted in Figure \ref{sample}. For each couple, we have access to the half-hour load consumption.

As we said previously, we want to cluster the relationship between two successive days.
In Figure \ref{uneSemaine}, we plot a week of load consumption. 

The regression model taken is $\boldsymbol{F}$ for the first day and $\boldsymbol{G}$ for the second day of each couple.
Besides, discretization of each day on $48$ points, every half-hour, is made available.
In our opinion, a linear model is not appropriate, as the behavior from the eve to the day depends on which day we consider: there is a difference between working days and weekend days, as involved in Figure \ref{uneSemaine}.

To apply our procedures, we project $\textbf{f}$ and $\textbf{g}$ onto wavelet basis. The symmlet-$4$ basis, at level $5$, is used.

We run our procedures with the number of clusters varying from $2$ to $6$.
Our both procedures select a model with $4$ components. The first one considers couples of weekdays, the second Friday-Saturday, the third component is Saturday-Sunday and the fourth
considers Sunday-Monday.
This result is faithful with the knowledge we have about these data. Indeed, working days have the same behavior, depending on the eve, whereas days off have not the same behavior, depending on working days and conversely.
Moreover, in the article of \cite{Misiti}, which also studied this example, they get the same classification.

\subsubsection{Tecator dataset}
This example deals with spectrometric data.
More precisely, a food sample has been considered, which contained finely chopped pure meat with different fat contents.
The data consist of a $100$-channel spectrum of absorbances in the wavelength range $850-1050$ nm and of the percentage of fat. We observe a sample of size $215$.
Those data have been studied in a lot of articles, refer for example to \cite{Ferraty}.
They work on different approaches. They test prediction and classification, supervised (where the fat content become a class, larger or smaller than $20 \%$), or not (ignoring the response variable). 
In this work, we focus on clustering data according to the relation between the fat content and the absorbance spectrum. 
We could not predict the response variable, because we do not know the class of a new observation. Estimate it is a difficult problem, in which we are not involved in this article.

We will take advantage of our procedures to know which coefficients, in the wavelet basis decomposition of the spectrum, are useful to describe the fat content.

The sample will be split into two subsamples, $165$ observations for the learning set and $50$ observations for the test set.
We split it to have the same marginal distribution for the response in each sample.

The spectrum is a function. We decompose it into the Haar basis, at level $6$.
Nevertheless, our model did not take into account a constant coefficient to describe the response.
Thereby, before running our procedure, we center the $y$ according to the learning sample and each function $\boldsymbol{x}_i$ for all observations in the whole sample.
Then, we could estimate the mean of the response by the mean $\hat{\mu}$ over the learning sample.

\begin{figure}[H]
\centering
 \includegraphics[scale=0.35,trim = 0cm 0cm 0cm 1.3cm, clip]{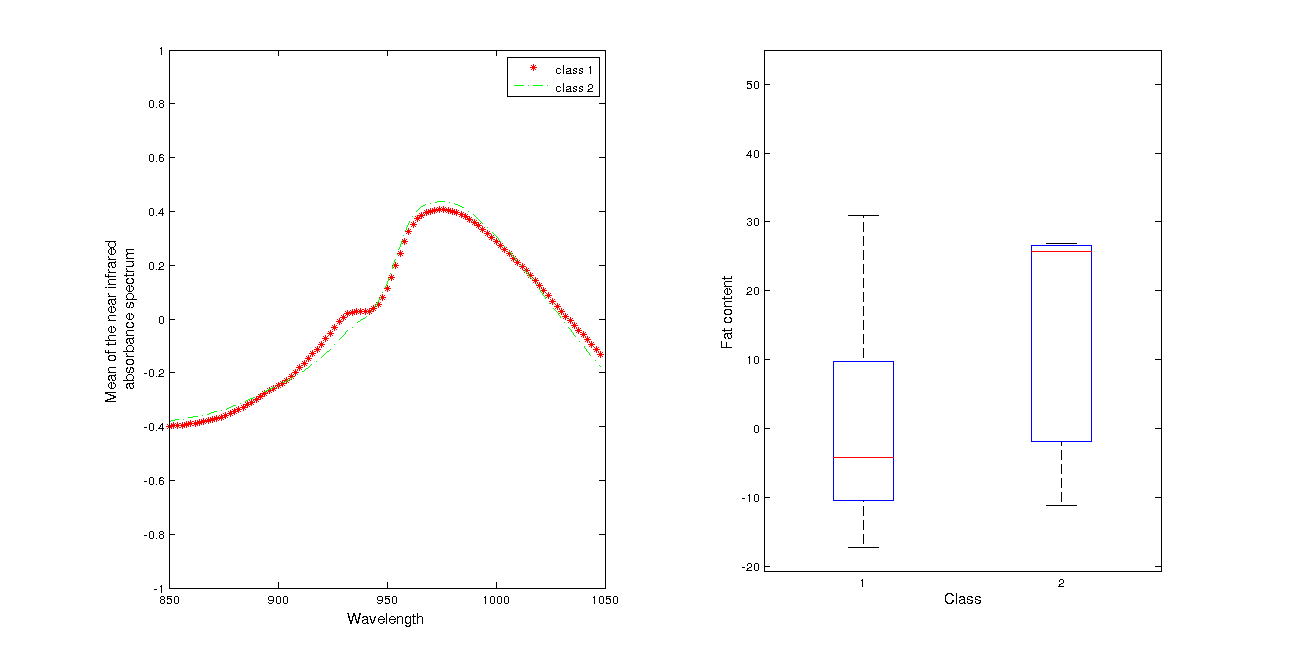}
 \caption[Summarized results for the model 1]{Summarized results for the model 1. The graph on the left is a candidate for representing each cluster, constructed by the mean of reconstructed spectrum over an a posteriori probability greater than $0.6$.
On the right side, we present the boxplot of the fat values in each class, for observations with an a posteriori probability greater than $0.6$. 
}
\label{Tecator1}
\end{figure}

\begin{figure}[H]
\centering
 \includegraphics[scale=0.35]{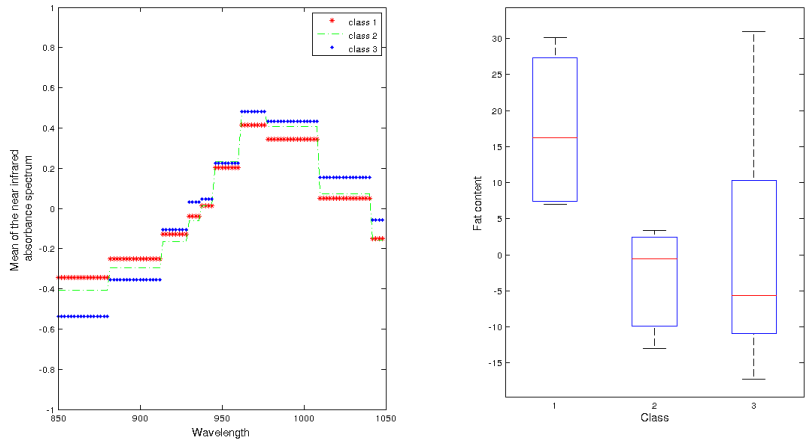}
 \caption[Summarized results for the model 1]{Summarized results for the model 2. The graph on the left is a candidate for representing each cluster, constructed by the mean of reconstructed spectrum over an a posteriori probability greater than $0.6$.
On the right side, we present the boxplot of the fat values in each class, for observations with an a posteriori probability greater than $0.6$. 
}
\label{Tecator2}
\end{figure}

We construct models on the training set by our procedure Lasso-MLE.
Thanks to the estimations, we have access to relevant variables and we could reconstruct signals by keeping only relevant variables.
We have also access to the a posteriori probability, which leads us to know which observation is,with high probability, in which cluster.
However, for some observations, the a posteriori probability do not ensure the clustering, because it is almost the same for different clusters.
The procedure selects two models that we describe here.
In Figures \ref{Tecator1} and \ref{Tecator2}, we represent clusters performed on the training set for the different models.
The graph on the left is a candidate for representing each cluster, constructed by the mean of spectrum over an a posteriori probability greater than $0.6$. We plot the reconstruction of the curve by keeping only the relevant variables in the wavelet decomposition.
On the right side, we present the boxplot of the fat values in each class, for observations with an a posteriori probability greater than $0.6$. 

The first model has two classes, which could be distinguish in the absorbance spectrum by the bump on wavelength around $940$ nm. 
The first cluster is dominating, with $\hat{\pi}_1=0.95$.
The fat content is smaller in the first cluster
than in the second cluster.
According to the signal reconstruction, we could see that almost all variables have been selected.
This model seems consistent according to the classification goal.

The second model has $3$ classes and we could remark different important wavelengths.
Around $940$ nm, there is a difference between classes, corresponding to the bump underline in the model 1.
 Around $970$ nm, there is also a difference between classes, with higher or smaller values.
The first class is dominating, with $\hat{\pi}_1=0.89$.
Just a few number of variables have been selected, which give to this model the understanding property of which coefficients are discriminating.


We could discuss about those models.
The first one selects only two classes, but almost all variables, whereas the second model has more classes and less variables: there is a trade-off between clusters and variable selection for the dimension reduction.
 
According to those classifications, we could compute the response according to the linear model.
We use two ways to compute $\hat{\boldsymbol{y}}$: either considering the linear model in the cluster selected by the MAP principle, or mixing estimations in each cluster thanks to these a posteriori probabilities.
We compute the Mean Absolute Percentage Error, $\text{MAPE}=\frac{1}{n} \sum_{i=1}^n |(\hat{\boldsymbol{y}}_i-\boldsymbol{y}_i) / \boldsymbol{y}_i|$.
Results are summarized in Table \ref{tableMPE}.
\begin{table}[H]
\centering
\begin{tabular}{|c|c|c|}
  \hline
   & Linear model in the class & Mixing estimation \\
   & with higher probability  & \\
  \hline
 Model 1 &  0.200   &  0.198\\
 Model 2 &  0.055    & 0.056 \\
  \hline
\end{tabular}
\caption{\label{tableMPE} Mean absolute percentage of error of the predicted value, for each model, for the learning sample.}
 \end{table}
Thus, we work on the test sample.
We use the response and the regressors to know the a posteriori of each observation. Then, using our models, we could compute the predicted fat values from the spectrometric curve, as before according to two ways, mixing or choosing the classes.
\begin{table}[H]
\centering
\begin{tabular}{|c|c|c|}
  \hline
   & Linear model in the class & Mixing estimation \\
   & with higher probability  & \\
  \hline
 Model 1 &  0.22196     &  0.21926\\
 Model 2 &  0.20492   & 0.20662 \\
  \hline
\end{tabular}
\caption{\label{tableMPE2} Mean absolute percentage of error of the predicted value, for each model, for the test sample.}
 \end{table}

 Because the models are constructed on the learning sample, the MAPE are lower than for the test sample. Nevertheless, results are similar,  saying that models are well constructed.
 This is particularly the case for the model 1, which is more consistent over a new sample.

To conclude this study, we could highlight the advantages of our procedure on these data. 
It provides a clustering of data, similar to the one done with supervised clustering in \cite{Ferraty}, but we could explain how this clustering is done.

This work has been done with the Lasso-MLE procedure. However, the same kind of results have been get with the Lasso-Rank procedure.
\section{Conclusion}
In this article, two procedures are proposed to cluster regression data.
Detecting the relevant clustering variables, they are especially designed for high-dimensional dataset.
We use an $\ell_1$-regularization procedure to select variables and then deduce a reasonable random model collection.
Thus, we recast estimations of parameters of these models into a general model selection problem.
These procedures are compared with usual criteria on simulated data: the BIC criterion used to select a model, the maximum-likelihood estimator and the oracle when we know it.
In addition, we compare our procedures to others on benchmark data. 

One main asset of those procedures is that it can be applied to functional dataset. We also develop this point of view.

 \section{Appendices}
 In those appendices, we develop computing for EM algorithm updating formulae in Section \ref{em}, for Lasso and maximum likelihood estimators and for low ranks estimators.
 In Section \ref{GroupLasso}, we extend our procedures with the Group-Lasso estimator to select relevant variables, rather than using the Lasso estimator.
 \subsection{EM algorithms} 
 \label{em}
 \subsubsection{EM algorithm for the Lasso estimator}

 Introduced by \cite{EM}, the EM (Expectation-Maximization) algorithm is used to compute maximum likelihood estimators, penalized or not.
 
 The expected complete negative log-likelihood is denoted by
$$Q(\theta|\theta')=-\frac{1}{n} E_{\theta'}(l_{c}(\theta,\textbf{X},\textbf{Y},\textbf{Z})|\textbf{X},\textbf{Y})$$
in which
\begin{align*}
l_{c} (\theta,\textbf{X}, \textbf{Y} ,\textbf{Z}) &= \sum_{i=1}^{n} \sum_{k=1}^{K} [\boldsymbol{Z}_i]_k \log\left( \frac{\det(\boldsymbol{P}_k)}{(2 \pi)^{q/2}} \exp \left(- \frac{1}{2} (\boldsymbol{P}_k \boldsymbol{Y}_i - \boldsymbol{X}_i\boldsymbol{\Phi}_k)^t (\boldsymbol{P}_k \boldsymbol{Y}_i - \boldsymbol{X}_i\boldsymbol{\Phi}_k) \right) \right)\\
&+ [\boldsymbol{Z}_i]_k \log(\pi_{k}) ;
\end{align*}

with $[\boldsymbol{Z}_i]_k$ are independent and identically distributed unobserved multinomial variables, showing the component-membership of the $i^{th}$ observation in the finite mixture regression model.

The expected complete penalized negative log-likelihood is
$$Q_{\text{pen}}(\theta|\theta')=Q(\theta|\theta')+ \lambda \sum_{k=1}^{K} \pi_{k} ||\boldsymbol{\Phi}_k||_1.$$
\paragraph{Calculus for updating formula}
\begin{itemize}
 \item E-step: compute  $Q(\theta|\theta^{(\text{ite})})$, or, equivalently, compute for $k \in \{1,\ldots,K\}$, $i \in \{1,\ldots,n\}$,
 \begin{align*}
\boldsymbol{\tau}^{\text{(ite)}}_{i,k}&=E_{\theta^{(\text{ite})}}([\boldsymbol{Z}_i]_k|\textbf{Y}) \\
&=\frac{\pi_{k}^{(\text{ite})} \det \boldsymbol{P}_k^{(\text{ite})} \exp \left( -\frac{1}{2}\left( \boldsymbol{P}_k^{(\text{ite})} \boldsymbol{Y}_i-\boldsymbol{X}_i \boldsymbol{\Phi}_k^{(\text{ite})}\right)^t\left(\boldsymbol{P}_k^{(\text{ite})} \boldsymbol{Y}_i-\boldsymbol{X}_i \boldsymbol{\Phi}_k^{(\text{ite})}\right) \right)}
{\sum_{r=1}^{K}\pi_{k}^{(\text{ite})}\det \boldsymbol{P}_k^{(\text{ite})} \exp \left( -\frac{1}{2}\left(\boldsymbol{P}_k^{(\text{ite})} \boldsymbol{Y}_i-\boldsymbol{X}_i \boldsymbol{\Phi}_k^{(\text{ite})}\right)^t\left(\boldsymbol{P}_k^{(\text{ite})} \boldsymbol{Y}_i-\boldsymbol{X}_i \boldsymbol{\Phi}_k^{(\text{ite})}\right) \right)}  
 \end{align*}
 This formula updates the clustering, thanks to the MAP principle.
 \item M-step: improve $Q_{\text{pen}}(\theta |\theta^{(\text{ite})})$.
 
 For this, rewrite the Karush-Kuhn-Tucker conditions. We have
 \begin{align}
 \label{AOptimiser}
&Q_\text{pen}(\theta|\theta^{\text{(ite)}}) \nonumber\\
=& - \frac{1}{n} \sum_{i=1}^{n} \sum_{k=1}^{K}  E_{\theta^{(\text{ite})}} \left( [\boldsymbol{Z}_i]_k \log \left( \frac{\det(\boldsymbol{P}_k)}{(2 \pi)^{q/2}} \exp \left( -\frac{1}{2} (\boldsymbol{P}_k \boldsymbol{Y}_i - \boldsymbol{X}_i \boldsymbol{\Phi}_k)^t(\boldsymbol{P}_k \boldsymbol{Y}_i - \boldsymbol{X}_i \boldsymbol{\Phi}_k) \right) \right) \left| \textbf{Y} \right. \right) \nonumber\\
&- \frac{1}{n} \sum_{i=1}^{n} \sum_{k=1}^{K} E_{\theta^{(\text{ite})}} \left( [\boldsymbol{Z}_i]_k \log \pi_{k} |\textbf{Y} \right) + \lambda \sum_{k=1}^{K} \pi_{k} ||\boldsymbol{\Phi}_k||_1  \nonumber\\
=& - \frac{1}{n} \sum_{i=1}^{n} \sum_{k=1}^{K} -\frac{1}{2} (\boldsymbol{P}_k \boldsymbol{Y}_i - \boldsymbol{X}_i \boldsymbol{\Phi}_k)^t (\boldsymbol{P}_k \boldsymbol{Y}_i - \boldsymbol{X}_i \boldsymbol{\Phi}_k)E_{\theta^{(\text{ite})}} \left( [\boldsymbol{Z}_i]_k |\textbf{Y} \right)\nonumber\\
&- \frac{1}{n} \sum_{i=1}^{n} \sum_{k=1}^{K}  \sum_{m=1}^q \log \left( \frac{[\boldsymbol{P}_k]_{m,m}}{\sqrt{2 \pi}} \right) E_{\theta^{(\text{ite})}} \left[ [\boldsymbol{Z}_i]_k |\textbf{Y} \right] \nonumber\\
&- \frac{1}{n} \sum_{i=1}^{n} \sum_{k=1}^{K} E_{\theta^{(\text{ite})}} \left( [\boldsymbol{Z}_i]_k  |\textbf{Y} \right)\log \pi_{k} + \lambda \sum_{k=1}^{K} \pi_{k} ||\boldsymbol{\Phi}_k||_1.
 \end{align}
Firstly, we optimize this formula with respect to $\boldsymbol{\pi}$: it is equivalent to optimize
$$-\frac{1}{n} \sum_{i=1}^{n} \sum_{k=1}^{K} \boldsymbol{\tau}_{i,k} \log(\pi_{k}) + \lambda \sum_{k=1}^{K} \pi_{k} ||\boldsymbol{\Phi}_k||_1.$$

We obtain
$$\pi_{k}^{(\text{ite}+1)}= \pi_{k}^{(\text{ite})}+t^{(\text{ite})} \left(\frac{\sum_{i=1}^{n} \boldsymbol{\tau}_{i,k}}{n} - \pi_{k}^{(\text{ite})}\right) ;$$
with $t^{(\text{ite})} \in (0,1]$ being the largest value in the grid $\{\delta^{l}, l\in \mathbb{N} \}$, with $0 < \delta <1$, such that the function is not increasing.

To optimize \eqref{AOptimiser} with respect to $(\boldsymbol{\Phi},\mathbf{P})$, we could rewrite the expression: it is similar to the optimization of
$$-\frac{1}{n} \sum_{i=1}^n \left(\boldsymbol{\tau}_{i,k} \sum_{m=1}^q \log([\boldsymbol{P}_k]_{m,m}) - \frac{1}{2} (\boldsymbol{P}_k \boldsymbol{\widetilde{Y}}_i - \boldsymbol{\widetilde{X}}_i \boldsymbol{\Phi}_k)^t(\boldsymbol{P}_k \boldsymbol{\widetilde{Y}}_i - \boldsymbol{\widetilde{X}}_i \boldsymbol{\Phi}_k) \right) + \lambda \pi_{k} ||\boldsymbol{\Phi}_k||_1$$
for all $k \in \{1,\ldots,K\}$, which is equivalent to the optimization of

$$-\frac{1}{n} n_k \sum_{m=1}^q \log( [\boldsymbol{P}_k]_{m,m}) + \frac{1}{2n} \sum_{i=1}^n \sum_{m=1}^q 
\left( [\boldsymbol{P}_k]_{m,m} [\boldsymbol{\widetilde{Y}}_i]_{k,m} - [\boldsymbol{\Phi}_k]_{m,.} [\boldsymbol{\widetilde{X}}_i]_{k,.} \right) ^2 + \lambda \pi_{k} ||\boldsymbol{\Phi}_k||_1 ;$$
where $n_k = \sum_{i=1}^n \boldsymbol{\tau}_{i,k}$.
The minimum in $[\boldsymbol{P}_k]_{m,m}$ is the function which vanishes its partial derivative with respect to $[\boldsymbol{P}_k]_{m,m}$:
$$- \frac{n_k}{n} \frac{1}{[\boldsymbol{P}_k]_{m,m}} + \frac{1}{2n}  \sum_{i=1}^n 2 [\boldsymbol{\widetilde{Y}}_i]_{k,m} \left( [\boldsymbol{P}_k]_{m,m} [\boldsymbol{\widetilde{Y}}_i]_{k,m} - [\boldsymbol{\Phi}_k]_{m,.} [\boldsymbol{\widetilde{X}}_i]_{k,.}\right)=0$$
for all $k \in \{1,\ldots,K\}$, for all $m \in \{1,\ldots,q\}$, which is equivalent to 

\begin{align*}
-1+\frac{1}{n_k} [\boldsymbol{P}_k]_{m,m}^2 \sum_{i=1}^n [\boldsymbol{\widetilde{Y}}_i]_{k,m}^2 - \frac{1}{n_k} [\boldsymbol{P}_k]_{m,m} \sum_{i=1}^n [\boldsymbol{\widetilde{Y}}_{i}]_{k,m} [\boldsymbol{\Phi}_k]_{m,.} [\boldsymbol{\widetilde{X}}_i]_{k,.}&=0\\
\Leftrightarrow -1 + [\boldsymbol{P}_k]_{m,m}^2 \frac{1}{n_k} || [\mathbf{\boldsymbol{\widetilde{Y}}}]_{k,m} ||_2^2 - [\boldsymbol{P}_k]_{m,m}\frac{1}{n_k} \langle[\widetilde{\textbf{Y}}]_{k,m}, [\boldsymbol{\Phi}_k]_{m,.} [\widetilde{\textbf{X}}]_{k,.} \rangle&=0.
\end{align*}
The discriminant is 
$$\Delta = \left( -\frac{1}{n_k} \langle [\widetilde{\textbf{Y}}]_{k,m},[\boldsymbol{\Phi}_k]_{m,.} [\widetilde{\textbf{X}}]_{k,.} \rangle \right)^2 - \frac{4}{n_k}  ||[\mathbf{\boldsymbol{\widetilde{Y}}}]_{k,m}||_2^2.$$
Then, for all $k \in \{1,\ldots,K\}$, for all $m \in \{1,\ldots,q\}$,
$$[\boldsymbol{P}_k]_{m,m}= \frac{n_k \langle [\widetilde{\textbf{Y}}]_{k,m}, [\boldsymbol{\Phi}_k]_{m,.} [\widetilde{\textbf{X}}]_{k,.} \rangle + \sqrt{\Delta}}{2 n_k ||[\mathbf{\boldsymbol{\widetilde{Y}}}]_{k,m}||_2^2}.$$

We could also look at the equation \eqref{AOptimiser} as a function of the variable $\boldsymbol{\Phi}$: according to the partial derivative with respect to 
$[\boldsymbol{\Phi}_k]_{m,j}$, we obtain for all $m \in \{1,\ldots,q\}$, for all $k \in \{1,\ldots,K\}$, for all $j \in \{1,\ldots,p\}$,
$$\sum_{i=1}^n [\boldsymbol{\widetilde{X}}_i]_{k,j} \left([\boldsymbol{P}_k]_{m,m} [\boldsymbol{\widetilde{Y}}_i]_{k,m} - \sum_{j_2 =1}^{p} [\boldsymbol{\widetilde{X}}_i]_{k,j_2} [\boldsymbol{\Phi}_k]_{m,j_2}\right) -n\lambda \pi_{k} \text{sgn}([\boldsymbol{\Phi}_k]_{m,j})=0,$$
where $\text{sgn}$ is the sign function.
Then, for all $k \in \{1,\ldots,K\}, j \in \{1,\ldots,p\}, m \in \{1,\ldots,q\}$,
$$[\boldsymbol{\Phi}_k]_{m,j} = \frac{\sum_{i=1}^{n} [\boldsymbol{\widetilde{X}}_i]_{k,j} [\boldsymbol{P}_k]_{m,m} [\boldsymbol{\widetilde{Y}}_i]_{k,m} - \sum_{\genfrac{}{}{0pt}{}{j_2=1}{j_2\neq j}}^{p} [\boldsymbol{\widetilde{X}}_i]_{k,j} [\boldsymbol{\widetilde{X}}_i]_{k,j_2} [\boldsymbol{\Phi}_k]_{m,j_2} - n\lambda \pi_{k} \text{sgn}([\boldsymbol{\Phi}_k]_{j,m})}{||[\mathbf{\boldsymbol{\widetilde{X}}}]_{k,j}||_{2}^2}.$$

Let, for all $k \in \{1,\ldots,K\}, j \in \{1,\ldots,p\}, m \in \{1,\ldots,q\}$, 
$$[\boldsymbol{S}_{k}]_{j,m}=-\sum_{i=1}^{n} [\boldsymbol{\widetilde{X}}_i]_{k,j} [\boldsymbol{P}_k]_{m,m} [\boldsymbol{\widetilde{Y}}_i]_{k,m} + \sum_{\genfrac{}{}{0pt}{}{j_2=1 }{ j_2\neq j}}^{p} [\boldsymbol{\widetilde{X}}_i]_{k,j} [\boldsymbol{\widetilde{X}}_i]_{k,j_2} [\boldsymbol{\Phi}_k]_{m,j_2}.$$
Then
\begin{align*}
[\boldsymbol{\Phi}_k]_{m,j}&= \frac{-[\boldsymbol{S}_{k}]_{j,m}- n\lambda \pi_{k} \text{sgn}([\boldsymbol{\Phi}_k]_{m,j})}{||[\mathbf{\boldsymbol{\widetilde{X}}}]_{k,j}||_{2}^2} \\
&= \left\{ \begin{array}{lll}
            &\frac{-[\boldsymbol{S}_{k}]_{j,m}+ n\lambda \pi_{k} }{||[\mathbf{\boldsymbol{\widetilde{X}}}]_{k,j}||_{2}^2} & \text{if }  [\boldsymbol{S}_{k}]_{j,m}>n \lambda \pi_{k}\\
            &-\frac{[\boldsymbol{S}_{k}]_{j,m}+ n\lambda \pi_{k}}{||[\mathbf{\boldsymbol{\widetilde{X}}}]_{k,j}||_{2}^2} & \text{if } [\boldsymbol{S}_{k}]_{j,m} < -n \lambda \pi_{k}\\
            &0 & \text{elsewhere.}
           \end{array}
\right.
\end{align*}
\end{itemize}
From these equalities, we could write the updating formulae.
 For $j \in \{1, \ldots, p\}, k \in \{1,\ldots,K \}$,  $m~\in~\{1,\ldots,q\}$, let
\begin{align*}
&[\boldsymbol{S}_{k}]_{j,m}^{(\text{ite})}=-\sum_{i=1}^{n} [\boldsymbol{\widetilde{X}}_i]_{k,j} [\boldsymbol{P}_k]^{(\text{ite})}_{m,m} [\boldsymbol{\widetilde{Y}}_i]_{k,m} + 
\sum_{\genfrac{}{}{0pt}{}{j_2=1}{ j_2\neq j}}^{p} [\boldsymbol{\widetilde{X}}_i]_{k,j} [\boldsymbol{\widetilde{X}}_i]_{k,j_2} [\boldsymbol{\Phi}_k]^{(\text{ite})}_{m,j_2} ;\\
&n_{k} = \sum_{i=1}^{n} \boldsymbol{\tau}_{i,k} ; \\
&([\boldsymbol{\widetilde{Y}}_{i}]_{k,.},[\boldsymbol{\widetilde{X}}_i]_{k,.}) = \sqrt{\boldsymbol{\tau}_{i,k}} (\boldsymbol{Y}_i,\boldsymbol{X}_i).
\end{align*}

\subsubsection{EM algorithm for the rank procedure}
To take into account the matrix structure, we want to make a dimension reduction on the rank of the mean matrix.
If we know to which cluster each sample belongs, we could compute the low rank estimator for linear model in each component.

Indeed, an estimator of fixed rank $r$ is known in the linear regression case: denoting $A^+$ the Moore-Penrose pseudo-inverse of $A$ and 
$[A]_r = U D_r V^t$ in which $D_r$ is obtained from $D$ by setting $(D_{r})_{i,i}=0$ for $i \geq r+1$, with $U D V^t$ the singular decomposition of $A$, 
if $\boldsymbol{Y}=\mathbf{B} X + \boldsymbol{\Sigma}$, an estimator of $\mathbf{B}$ with rank $r$
is $\hat{\mathbf{B}}_r = [(\mathbf{x}^t \mathbf{x})^{+} \mathbf{x}^t \mathbf{y}]_r$.

 We do not know the clustering of the sample, but the E-step in the EM algorithm computes it.
 
 \label{EM2}
We suppose in this case that $\boldsymbol{\Sigma}_k$ and $\pi_k$ are known, for all $k \in \{1,\ldots,K\}$.
We use this algorithm to determine $\boldsymbol{\Phi}_k$, for all $k \in \{1,\ldots,K\}$, with ranks fixed to $\textbf{R} = (R_1, \ldots, R_K)$.
\begin{itemize}
 \item E-step: compute for $k \in \{1,\ldots,K\}$, $i \in \{1,\ldots,n\}$,
 \begin{align*}
\boldsymbol{\tau}_{i,k}&=E_{\theta^{(\text{ite})}}([\boldsymbol{Z}_i]_k|Y) \\
&=\frac{\pi_{k}^{(\text{ite})} \det \boldsymbol{P}_k^{(\text{ite})} \exp \left( -\frac{1}{2}\left( \boldsymbol{P}_k^{(\text{ite})} \boldsymbol{y}_i-\boldsymbol{x}_i \boldsymbol{\Phi}_k^{(\text{ite})}\right)^t\left(\boldsymbol{P}_k^{(\text{ite})} \boldsymbol{y}_i-\boldsymbol{x}_i \boldsymbol{\Phi}_k^{(\text{ite})}\right) 
\right)}{\sum_{r=1}^{K}\pi_{k}^{(\text{ite})}\det \boldsymbol{P}_k^{(\text{ite})} \exp \left( -\frac{1}{2}\left(\boldsymbol{P}_k^{(\text{ite})} \boldsymbol{y}_i-\boldsymbol{x}_i \boldsymbol{\Phi}_k^{(\text{ite})}\right)^t
\left(\boldsymbol{P}_k^{(\text{ite})} \boldsymbol{y}_i-\boldsymbol{x}_i \boldsymbol{\Phi}_k^{(\text{ite})}\right) \right)}  
 \end{align*}
 \item M-step: assign each observation to its estimated cluster, by the MAP principle applied thanks to the E-step. We say that $\boldsymbol{Y}_i$ comes from component number $\underset{k \in \{ 1,\ldots, K\}}{ \operatorname{argmax}}\boldsymbol{\tau}^{\text{(ite)}}_{i,k}$.
 Then, we can define $\widetilde{\mathbf{B}_k}^{\text{(ite)}} = (\mathbf{x}^t_{|k} \mathbf{x}_{|k})^{-1} \mathbf{x}^t_{|k} \mathbf{y}_{|k}$, in which $\mathbf{x}_{|k}$ and $\mathbf{y}_{|k}$ 
 correspond to the observations belonging to the cluster $k$.
 We decompose $\widetilde{\mathbf{B}}_{k}^{(\text{ite})}$ in singular value: $\widetilde{\mathbf{B}}_{k}^{(\text{ite})} = U S V^t$.
 Then, the estimator is $\hat{\mathbf{B}}_{k}^{(\text{ite})} = U \boldsymbol{S}_{R(k)} V^t$.
\end{itemize}

\subsection{Group-Lasso MLE and Group-Lasso Rank procedures}
\label{GroupLasso}
One way to perform those procedures is to consider the Group-Lasso estimator rather than the Lasso estimator to select relevant variables.
Indeed, this estimator is more natural, according to the relevant variable definition.
Nevertheless, results are very similar, because we select grouped variables in both case, selected by the Lasso or by the Group-Lasso estimator.
In this section, we describe our procedures with the Group-Lasso estimator, which could be understood as an improvement of our procedures.
\subsubsection{Context - definitions}
Our both procedures take advantage of the Lasso estimator to select relevant variables and to reduce the dimension in case of high-dimensional dataset.
First, recall what is a relevant variable.

\begin{definition}
A variable indexed by $(m,j) \in \{1,\ldots,q\} \times \{1,\ldots,p\}$ is \emph{irrelevant} for the clustering if 
$$[\boldsymbol{\Phi}_{1}]_{m,j} = \ldots = [\boldsymbol{\Phi}_{K}]_{m,j}=0.$$
A \emph{relevant} variable is a variable which is not irrelevant.
We denote by $J$ the set of relevant variables.
\end{definition}

According to this definition, we could introduce the Group-Lasso estimator.

\begin{definition}
The \emph{Lasso estimator} for mixture regression models with regularization parameter $\lambda \geq 0$ is defined by 
\begin{equation*}
 \hat{\theta}^{\text{Lasso}}(\lambda) := \underset{\theta \in \Theta_K}{\operatorname{argmin}} \left\{ -\frac{1}{n} l_{\lambda}(\theta) \right\} ;
 \label{lassomultidim2}
 \end{equation*}
where
$$- \frac{1}{n} l_\lambda (\theta) = -\frac{1}{n} l(\theta) + \lambda \sum_{k=1}^{K} \pi_k ||\boldsymbol{\Phi}_k||_1 ;$$
where $||\boldsymbol{\Phi}_k||_1 = \sum_{j=1}^p \sum_{m=1}^q |[\boldsymbol{\Phi}_k]_{m,j}|$ and with $\lambda$ to specify.
\end{definition}
It is the estimator used in the both procedures described in previous parts.
\begin{definition}
The \emph{Group-Lasso estimator} for mixture regression models with regularization parameter $\lambda \geq 0$ is defined by 
\begin{equation*}
 \hat{\theta}^{\text{Group-Lasso}}(\lambda) := \underset{\theta \in \Theta_K}{\operatorname{argmin}} \left\{ -\frac{1}{n} \widetilde{l}_{\lambda}(\theta) \right\} ;
 \label{GroupLassoeq}
 \end{equation*}
where
$$- \frac{1}{n} \widetilde{l}_\lambda (\theta) = -\frac{1}{n} l(\theta) + \lambda \sum_{j=1}^{p} \sum_{m=1}^q \sqrt{k}||[\boldsymbol{\Phi}]_{m,j}||_2 ;$$
where $||[\boldsymbol{\Phi}]_{m,j}||_2^2 = \sum_{k=1}^K |[\boldsymbol{\Phi}_k]_{m,j}|^2$ and with $\lambda$ to specify.
\end{definition}
This Group-Lasso estimator has the advantage to shrink to zero grouped variables rather than variables one by one.
It is consistent with the relevant variable definition.

However, depending on the dataset, it could be interesting to look for which variables are equal to zero first.
One way could be to extend this work with Lasso-Group-Lasso estimator, described for the example for the linear model in \cite{SparseGroupLasso}.

Let us describe two additional procedures, which will use the Group-Lasso estimator rather than the Lasso estimator to detect relevant variables.

 \subsubsection{Group-Lasso-MLE procedure}
This procedure is decomposed into three main steps: we construct a model collection, then in each model we compute the maximum likelihood estimator and we select the best one among all the models.

The first step consists of constructing a collection of models $\{\mathcal{H}_{(K,\widetilde{J})}\}_{(K,\widetilde{J}) \in \mathcal{M}}$ 
in which $\mathcal{H}_{(K,\widetilde{J})}$ is defined by 
\begin{equation}
\label{modele h bis}
\mathcal{H}_{(K,\widetilde{J})} = \left\{  y \in \mathbb{R}^q | x \in \mathbb{R}^p \mapsto h_{\theta}(y|x) \right\} ;
 \end{equation}
where 
$$  h_{\theta}(y|x)= \sum_{k=1}^{K} \frac{\pi_{k} \det(\boldsymbol{P}_k)}{(2 \pi)^{q/2}} \exp \left( -\frac{(\boldsymbol{P}_ky-\boldsymbol{\Phi}_k^{[\widetilde{J}]} x )^t(\boldsymbol{P}_ky-\boldsymbol{\Phi}_k^{[\widetilde{J}]}  x)}{2} \right),$$
and 
$$\theta=(\pi_1,\ldots, \pi_K,\boldsymbol{\Phi}_1,\ldots, \boldsymbol{\Phi}_K, \boldsymbol{P}_1,\ldots,\boldsymbol{P}_K) \in \Pi_K \times \left( \mathbb{R}^{q\times p} \right)^K \times \left(\mathbb{R}_+^{q} \right)^K.$$

The model collection is indexed by $ \mathcal{M}= \mathcal{K} \times \widetilde{\mathcal{J}}$. Denote by $\mathcal{K} \subset \mathbb{N}^*$ the possible number of components.
We could bound $\mathcal{K}$ without loss of estimation. Denote also $\widetilde{\mathcal{J}}$ a collection of subsets of $\{1,\ldots,q\} \times \{1,\ldots,p\}$, constructed by the Group-Lasso estimator.

To detect the relevant variables and construct the set $\widetilde{J} \in \widetilde{\mathcal{J}}$, we will use the Group-Lasso estimator defined by \eqref{GroupLassoeq}.
In the $\ell_1$-procedures, the choice of the regularization parameters is often difficult: fixing the number of components $K \in \mathcal{K}$,
we propose to construct a data-driven grid  $G_K$ of regularization parameters by using the updating formulae of the mixture parameters in the EM algorithm.

Then, for each $\lambda \in G_K$, we could compute the Group-Lasso estimator defined by
$$\hat{\theta}^{\text{Group-Lasso}} = \underset{\theta \in \Theta_K}{ \operatorname{argmin}} 
\left\{ -\frac{1}{n} \sum_{i=1}^n \log (h_\theta(\boldsymbol{y}_i|\boldsymbol{x}_i)) + \lambda \sum_{j=1}^{p} \sum_{m=1}^q \sqrt{K}||[\boldsymbol{\Phi}]_{m,j}||_2 \right\} .$$ 
For a fixed number of mixture components $K \in \mathcal{K}$ and a regularization parameter $\lambda$, we could use a generalized EM algorithm to approximate this estimator.
Then, for each $K \in \mathcal{K}$ and for each $\lambda \in G_K$, we have constructed the set of relevant variables $\widetilde{J}_\lambda$.
We denote by $\widetilde{\mathcal{J}}$ the collection of all these sets.

The second step consists of approximating the MLE 
$$\hat{h}^{(K,\widetilde{J})}= \underset{t \in \mathcal{H}_{(K,\widetilde{J})}}{ \operatorname{argmin}} \left\{ -\frac{1}{n} \sum_{i=1}^n \log (t(\boldsymbol{y}_i|\boldsymbol{x}_i)) \right\} ; $$ 
using the EM algorithm for each model $(K,\widetilde{J})\in \mathcal{M}$. 

The third step is devoted to model selection. We use the slope heuristic described in \cite{BirgeMassart}. 

\subsubsection{Group-Lasso-Rank procedure}
We propose a second procedure to take into account the matrix structure.
For each model belonging to the collection $\mathcal{H}_{(K,\widetilde{J})}$, a subcollection is constructed, varying the rank of $\boldsymbol{\Phi}$.
Let us describe this procedure.

As in the Group-Lasso-MLE procedure, we first construct a collection of models, thanks to the $\ell_1$-approach.
We obtain an estimator for $\theta$, denoted by $\hat{\theta}^{\text{Group-Lasso}}$, for each model belonging to the collection.
We could deduce the set of relevant variables, denoted by $\widetilde{J}$ and this for all $K \in \mathcal{K}$: we deduce $\widetilde{\mathcal{J}}$ the collection of set of relevant variables.

The second step consists in constructing a subcollection of models with rank sparsity, denoted by 
$$\{\check{\mathcal{H}}_{(K,\widetilde{J},R)}\}_{(K,\widetilde{J},R) \in \widetilde{\mathcal{M}}}.$$ 
The model $\{\check{\mathcal{H}}_{(K,\widetilde{J},R)}\}$ has $K$ components, the set $\widetilde{J}$ for active variables and $R$ is the vector
of the ranks of the matrix of regression coefficients in each group:

\begin{equation}
\label{modele htilde bis}
\check{\mathcal{H}}_{(K,\widetilde{J},R)}= \left\{ \boldsymbol{y} \in \mathbb{R}^q | \boldsymbol{x} \in \mathbb{R}^p \mapsto h^{(K,\widetilde{J},R)}_{\theta}(\boldsymbol{y}|\boldsymbol{x})   \right\}
 \end{equation}
where
\begin{align*}
h^{(K,\widetilde{J},R)}_{\theta}(\boldsymbol{y}|\boldsymbol{x}) =& \sum_{k=1}^{K} \frac{\pi_{k} \det(\boldsymbol{P}_k)}{(2 \pi)^{q/2}} \exp \left( -\frac{(\boldsymbol{P}_k\boldsymbol{y}-(\boldsymbol{\Phi}_k^{R_k})^{[\widetilde{J}]} \boldsymbol{x} )^t(\boldsymbol{P}_k\boldsymbol{y}-(\boldsymbol{\Phi}_k^{R_k})^{[\widetilde{J}]} \boldsymbol{x} )}{2} \right) ;\\
\theta=&(\pi_1,\ldots, \pi_K,\boldsymbol{\Phi}_1^{R_1},\ldots, \boldsymbol{\Phi}_K^{R_K}, \boldsymbol{P}_1,\ldots,\boldsymbol{P}_K) \in \Pi_K \times \Psi_K^R \times \left(\mathbb{R}_+^{q} \right)^K ;\\
\Psi_K^R=& \left\{(\boldsymbol{\Phi}_1^{R_1},\ldots, \boldsymbol{\Phi}_K^{R_K}) \in \left(\mathbb{R}^{q\times p} \right)^K | \text{Rank}(\boldsymbol{\Phi}_1) = R_1,\ldots, \text{Rank}(\boldsymbol{\Phi}_K)=R_K \right\} ; 
\end{align*}

and $\widetilde{\mathcal{M}}^R = \mathcal{K} \times \widetilde{\mathcal{J}} \times \mathcal{R}$. 
Denote by $\mathcal{K}\subset \mathbb{N}^*$ the possible number of components, $\widetilde{\mathcal{J}}$ a collection of subsets of $\{1,\ldots,q\} \times \{1,\ldots,p\}$ and $\mathcal{R}$ the set of vectors of size $K \in \mathcal{K}$ with rank values for each mean matrix.
We could compute the MLE under the rank constraint thanks to an EM algorithm.
Indeed, we could constrain the estimation of $\boldsymbol{\Phi}_k$, for all $k$, to have a rank equal to $R_k$, in keeping only the $R_k$ largest singular values.
More details are given in Section \ref{EM2}.
It leads to an estimator of the mean with row sparsity and low rank for each model.

\end{document}